\documentclass[12pt]{amsart}
\usepackage{amsthm,amsfonts,amssymb,amscd,mathrsfs}
\usepackage{bezier,longtable,amstext}
\usepackage{amssymb}

\usepackage[all]{xy}

\newcommand{\EXT}{{\mathcal Ext}}

\newcommand\codim{\text{codim}}

\newcommand{\const}{\operatorname{const}\nolimits}

\newcommand\Ext{\operatorname{Ext}\nolimits}

\newcommand{\Hilb}{\operatorname{Hilb}\nolimits}

\newcommand\Hom{\operatorname{Hom}\nolimits}
\newcommand\id{\text{id}}

\newcommand\rk{\text{rk}}

\newcommand\Span{\text{Span}}

\newcommand{\PP}{\mathbb{P}}

\newcommand{\GG}{\mathbf{G}}

\renewcommand{\SS}{\mathcal{S}}
\newcommand{\OO}{\mathcal{O}}
\newcommand{\FF}{\mathcal{F}}
\renewcommand{\phi}{\varphi}

\newcommand {\ZZ} {\mathbf {Z}}
\newcommand {\CC} {\mathbb {C}}

\newcommand{\refle}[1]{Lemma \ref{#1}}
\newcommand{\refth}[1]{Theorem \ref{#1}}

\newcommand{\refsec}[1]{Section~ \ref{#1}}
\newcommand{\refeq}[1]{(\ref{#1})}

\newtheorem{theorem}{Theorem}[section]
\newtheorem{proposition}[theorem]{Proposition}

\newtheorem{lemma}[theorem]{Lemma}
\newtheorem{corollary}[theorem]{Corollary}

\theoremstyle{definition}

\newtheorem{proposition-definition}[theorem]{Proposition-Definition}

\begin{document}

\title{Triviality of vector bundles on sufficiently twisted ind-Grassmannians}

\author[I.Penkov]{\;Ivan~Penkov}

\address{
Jacobs University Bremen \\
School of Engineering and Science,
Campus Ring 1,
28759 Bremen, Germany}
\email{ivanpenkov@jacobs-university.de}

\author[Tikhomirov]{\;Alexander S.~Tikhomirov}

\address{
Department of Mathematics\\
State Pedagogical University\\
Respublikanskaya Str. 108
\newline 150 000 Yaroslavl, Russia}
\email{astikhomirov@mail.ru}


\maketitle

\thispagestyle{empty}

\begin{abstract}
Twisted ind-Grassmannians are ind-varieties $\GG$ obtained as direct limits of Grassmannians $G(r_m,V^{r_m})$, for $m\in\ZZ_{>0}$, under embeddings $\phi_m:G(r_m,V^{r_m})\to G(r_{m+1}, V^{r_{m+1}})$ of degree greater than one. It has been conjectured in \cite{PT} and \cite{DP} that any vector bundle of finite rank on a twisted ind-Grassmannian is trivial. We prove this conjecture under the assumption that the ind-Grassmannian $\GG$ is sufficiently twisted, i.e. that $\lim_{m\to\infty}\frac{r_m}{\deg \phi_1\dots\deg\phi_m}=0$. 

2000 Mathematics Subject Classification, Primary 14M15, (Secondary 14J60, 32L05).
\end{abstract}
\section{Introduction}\label{sec1}
\label{Introduction}
{\it Ind-Grassmannians} are ind-varieties defined by chains of embeddings
\begin{equation}\label{eq1}
G(r_1,V^{n_1})\stackrel{\phi_1}{\hookrightarrow}G(r_2,V^{n_2})
\stackrel{\phi_2}{\hookrightarrow}\dots\stackrel{\phi_{m-1}}{\hookrightarrow}G(r_m,V^{n_m})
\stackrel{\phi_m}{\hookrightarrow}\dots,
\end{equation}
where $G(r,V)$ denotes the Grassmanian of $r$-planes in a finite dimensional vector space $V$. Any of the embeddings $\phi_m$ has a well defined degree $\deg \phi_m$, and the ind-Grassmannian defined by \refeq{eq1} is {\it twisted} if $\deg\phi_m>m$ for infinitely many
indices. In the special case when $r_m=1$ and $\deg\phi_m=1$ for all $m$,  the study of finite rank vector bundles on ind-Grassmannians goes
back to W.Barth, A.Van de Ven and A.N.Tyurin, \cite{BV}, \cite{T}. In this case the ind-Grassmannian is just the
infinite projective space $\mathbb{P}^\infty$, and the remarkable Barth-Van de Ven-Tyurin
Theorem claims that any vector bundle of finite rank on $\mathbb{P}^\infty$ is isomorphic to a
direct sum of line bundles. Historically, this is the first manifestation of a general phenomenon that seems
to take place for ind-varieties defined via sequences of embeddings similar to
(\ref{eq1}) with $G(r_m,V_m)$ replaced by arbitrary compact homogeneous spaces: in all cases known, the
restriction of any finite rank vector bundle on the ind-variety to a large enough finite dimensional homogeneous subspace is equivariant. Around the same time this phenomenon occured also in the important work of E. Sato who gave an independent proof of the Barth-Van de Ven-Tyurin Theorem, \cite{S1}. Shortly after that Sato established a more general result which applies in particular to the ind-Grassmannian $\mathbf{G}(r,V)$ of $r$-planes in a countable dimensional vector space $V$, \cite{S2}.

More recently the subject has been revisited in the papers \cite{DP}, \cite{CT} and \cite{PT}. In particular, in \cite{PT} a general conjecture about finite rank vector bundles on
ind-Grassmannians has been stated. In fact, as we show in \cite{PT}, if the ind-Grassmannian is not twisted (which is easily seen to be equivalent to assuming that $\deg\phi_m=1$ for all $m$), this conjecture is a relatively straightforward corollary of Sato's result. This leaves open the case of a twisted ind-Grassmannian, in which case the conjecture claims simply that finite rank vector bundle on such an ind-Grassmannian is trivial. So far this latter conjecture is established in the following three cases: for a rank two bundle on any twisted ind-Grassmannian \cite{PT}, for any finite rank bundle on any twisted projective ind-space (a twisted projective ind-space can be defined via the sequence (\ref{eq1}) where $r_m=1$ and $\deg\phi_m>1$ for all $m$) \cite{DP} and for an arbitrary finite rank bundle on some special twisted ind-Grassmannians (here $\phi_m$ are twisted extensions as defined in \cite{DP}).

In the present paper we consider the case of arbitrary finite rank vector bundle on arbirary twisted ind-Grassmannians satisfying the condition $\dim r_m=\const$ for all $m$. In fact, we work with a more general class of twisted ind-Grassmannians which we call {\it sufficiently twisted}. They are defined via the condition
\begin{equation}\label{tame}
\lim_{m\to\infty}\frac{r_m}{\deg\phi_1...\deg\phi_{m-1}}=0.
\end{equation}
Our main idea is that the relatively simple proof of the conjecture in the case of a twisted ind-projective space, \cite{DP}, admits an interesting generalization. More precisely, the original method is based on the study of certain morphisms of $\PP^1\times\PP^1$ into larger and larger projective spaces. Under the assumption that a twisted infinite projective space admits a non-trivial vector bundle, one pulls it back to $\PP^1\times\PP^1$, and for a sufficiently large projective space the pull-back is forced to have numerical invariants which yield a contradiction. The main technical achievement of the present paper is the introduction of an appropriate class of maps of $\mathbb{P}^1\times\mathbb{P}^1$ into a twisted ind-Grassmannian (defined in terms of a construction of a certain vector bundles on $\mathbb{P}^1\times\mathbb{P}^1$, see Section \ref{vector bdl} below) and the corresponding extension of certain estimates in \cite{DP} concerning lines in $\mathbb{P}^n$ to Segre curves of degree $2r$ in $G(r,V)$. The only limitation of this more general method has to do with fact that it is still based on the specific properties of the surface $\mathbb{P}^1\times\mathbb{P}^1$. This explains the condition (\ref{tame}).

\textbf{Acknowledgement. }We acknowledge the support and hospitality of the Max Planck Institute
for Mathematics in Bonn where the present paper was conceived. A. S. T. also acknowledges
partial support from Jacobs University Bremen.

\section{Notation and Conventions}\label{sec2}
Our notation is mostly standard. The ground field is $\CC$. All vector bundles considered are
assumed to have finite rank. We do not make a distinction between locally free sheaves of finite rank and vector bundles. If $\FF$ is a sheaf of $\OO_X$-modules on an algebraic variety $X$, $\FF^n$
denotes the direct sum of $n$ copies of $\FF$, $H^i(\FF)$ denotes the $i^{th}$ cohomology group of
$\FF$, $h^i(\FF):=\dim H^i(\FF)$, and $\FF^\vee$ stands for the dual bundle, i. e.
$\FF^\vee:= \mathcal{H}om_{\OO_X}(\FF,\OO_X)$. If $Z\subset X$ is a subvariety, $I_{Z,X}$ denotes
the sheaf of ideals corresponding to $Z$.

By $G(r,V)$ we denote the Grassmannian of $r$-dimensional subspaces of a vector space $V$;
unless the contrary is stated explicitly, we assume that $\dim V<\infty$,
$r\neq 1, r\neq\dim V-1$.

By a \textit{rational curve} we always mean a curve isomorphic to $\PP^1$, i. e.
for convenience we assume a rational curve to be smooth. If $C$ is a rational curve,
$\OO_C(i)$ stands for a line bundle on $C$ with first Chern class equal to $i\in\ZZ$.
A \emph{line} in $G(r,V)$ is a rational curve of degree 1 and is determined by a flag of
$V_1\subset V_2$ of subspaces in $V$ with $\dim V_1=r-1$, $\dim V_2=r+1$.

If $C\subset X$ is a rational curve in an algebraic variety $X$ and $E$ is a vector bundle on
$X$, then by a classical theorem of Grothendieck, $\displaystyle E_{|C}$ is isomorphic to
$\bigoplus_i\OO_C(d_i)$ for some $d_1\geq d_2\geq\dots\geq d_{\rk E}$. We call the ordered
$\rk E$-tuple $(d_1,\dots,d_{\rk E})$ \emph{the splitting type} of $E_{|C}$ and denote it by
$\mathbf{d}_E(C)$.

We call a curve $C=\overset{r}{\underset{i=1}\cup} C_i$, where $C_i$ are rational curves,
a \textit{chain of rational curves}, if, for each $i<r$, the intersection $C_i\cap C_{i+1}$
is a transversal intersection at a single point and there are no other intersections
of the curves $C_i$.
If $C$ is a chain of rational curves,  $\mathcal{O}_C(n_1,...,n_r)$ denotes a line bundle on $C$ such that
$\mathcal{O}_C(n_1,...,n_r)_{|C_i}\simeq\mathcal{O}_{C_i}(n_i)$.

Finally, under a \emph{partition} of $n\in\ZZ_{>0}$ (respectively, a \emph{strict partition}
of $n$) we understand a $k$-tuple $(n_1,\dots,n_k)\in \ZZ^k_{\geq 0}$ (respectively,
$(n_1,\dots,n_k)\in \ZZ^k_{> 0}$) with $\sum_{i=1}^k n_i=n$.

\vspace{0.5cm}

\section{An estimate for $\mathbf{D}_E(C)$}

\vspace{0.5cm}

For a vector bundle $E_C$ on a rational curve $C$ with splitting type
$\mathbf{d}(E_C)=(\mathbf{d}_1(E_C),\dots,\mathbf{d}_{\rk E}(E_C))$, set
$\mathbf{D}(E_C):=\mathbf{d}_1(E_C)-\mathbf{d}_{\rk E}(E_C)$.
Our objective in this section is to prove the following theorem.

\begin{theorem}\label{theorem 5.1}
Let $\pi:X\to B$ be a flat family whose fibers are chains of rational curves and whose generic fiber is a rational curve. Assume that, for a point $0\in B$, $C_0:=\pi^{-1}(0)$ is a chain of rational curves $C_1\cup...\cup C_r,\ r\ge1$. Let $E$ be a vector bundle on $X$. Then there exists a neighbourhood $U$ of  the point $0$ in $B$ such that for any $t\in U$ for which $C_t=\pi^{-1}(t)$ is a rational curve, one has

\begin{equation}\label{ineq4}
\mathbf{D}(E_{|C_t})\le\overset{r}{\underset{i=1}\sum}\mathbf{D}(E_{|C_i}).
\end{equation}
\end{theorem}
\begin{proof}
The proof is based on two auxiliary results, namely Corollary \ref{exists} and Lemma \ref{ineq2},
which we prove later on in this section. First, Corollary \ref{exists} implies that for any
integers $n_1,...,n_r$ there exists a neighbourhood $U'\subset B$ of the point $0$ and a line
bundle $\mathcal{L}$ on $\pi^{-1}(U')$ such that
$\mathcal{L}_{|C_0}\simeq\mathcal{O}_{C_0}(n_1,...,n_r)$. Hence
$E_{|C_0}\otimes \mathcal{L} \simeq E_{|C_0}(n_1,...,n_r)$ and
$E_{|C_t}\otimes \mathcal{L} \simeq E_{|C_t}(n_1+...+n_r)$ for any $t\in U'$ for which
$C_t=\pi^{-1}(t)$ is a rational curve. By semicontinuity,
\begin{equation}\label{ineq50}
h^0((E_{|C_0})(n_1,...,n_r))\ge h^0((E_{|C_t})(n_1+...+n_r)).
\end{equation}
Therefore, for $n_i=-\mathbf{d}_1(E_{|C_i})-\delta_i$, where $\delta_i$ are as in Lemma
 \ref{ineq2}, the inequality (\ref{ineq50}) and Lemma \ref{ineq2} imply
\begin{equation}\label{ineq500}
h^0(E_{|C_t}(-\overset{r}{\underset{i=1}\sum}\mathbf{d}_1(E_{|C_i})
-\overset{r}{\underset{i=1}\sum}\delta_i))=0.
\end{equation}
In particular, (\ref{ineq500}) holds for the following $r$ choices of  $\delta_1,...,\delta_r$:
$\delta_{i_0}=1,...,\delta_j=0$ for $j\ne i_0$, $i_0$ running from $1$ to $r$.
Therefore, for $t$ in the intersection $U$ of the corresponding $r$ neighbourhoods $U'$
we have
\begin{equation}\label{ineq5}
\mathbf{d}_1(E_{|C_t})\le\overset{r}{\underset{i=1}\sum}\mathbf{d}_1(E_{|C_i}).
\end{equation}
Since
$\mathbf{d}_{\rk E}(E)=-\mathbf{d}_1(E^\vee)$, inequality (\ref{ineq5}) applied
to $E^\vee$ instead of $E$ yields
\begin{equation}\label{ineq51}
\mathbf{d}_{\rk E}(E_{|C_t})\ge\overset{r}{\underset{i=1}\sum}\mathbf{d}_{\rk E}(E_{|C_i})
\end{equation}
for $t\in U$. The desired inequality (\ref{ineq4}) follows from (\ref{ineq5}) and
(\ref{ineq51}).
\end{proof}

We now proceed to the auxiliary statements used above. Given a strict partition
$(n_1,\dots,n_r)$ of $n\in \ZZ_{>0}$, we define a \emph{polarized chain} (of rational cruves)
as a pair $(C,\OO_C(n_1,\dots,n_r))$, where $C=C_1\cup...\cup C_r$ is a chain of rational
curves.

\begin{lemma}\label{embed}
For any strict partition $(n_1,...,n_r)$ of $n\in\ZZ_{>0}$ and any polarized chain
$(C,\ \mathcal{O}_C(n_1,...,n_r))$ there exists a linearly normal \footnote{Recall that $i$
is linearly normal if $\Span(i(C))=\mathbb{P}^n$.} embedding $i:C\hookrightarrow\mathbb{P}^n$
such that $\mathcal{O}_C(n_1,...,n_r))\simeq i^*\mathcal{O}_{\mathbb{P}^n}(1).$
\end{lemma}

\begin{proof}
We use induction on $r$. For $r=1$ the desired embedding $i:C\hookrightarrow\mathbb{P}^n$ is
clearly given by the complete linear series $|\mathcal{O}_{\mathbb{P}^1}(n_1)|$.
Assume now that the claim is true for $r-1$. If we decompose $C$ as $C=C'\cup C_r$,
where $C':=C_1\cup...\cup C_{r-1}$, and set $n':=n-n_r$, then by the induction assumption
there exists a linearly normal embedding
$i':C'\hookrightarrow\mathbb{P}^{n'}$ such that $\mathcal{O}_C(n_1,...,n_{r-1}))=
i'^*\mathcal{O}_{\mathbb{P}^{n'}}(1).$
Next, consider the linearly normal embedding $i_r:C_r\hookrightarrow\mathbb{P}^r$ by the
complete linear series $|\mathcal{O}_{\mathbb{P}^1}(n_r)|$ and embed the spaces
$\mathbb{P}^{n'}=\Span(i'(C'))$ and $\mathbb{P}^r=\Span(i_r(C_r))$ into the projective space
$\mathbb{P}^n$ in such a way that their intersection $\mathbb{P}^{n'}\cap\mathbb{P}^r$ in $\mathbb{P}^n$
is a point. We may assume, after possible projective linear transformations of
$\mathbb{P}^{n'}$ and $\mathbb{P}^r,$ that this point equals $i'(C')\cap i_r(C_r)$.
Thus we obtain an embedding $i:C\hookrightarrow\mathbb{P}^n$ such that
$i_{|C'}=i'$, $i_{|C_r}=i_r$ and, by the construction,
$\mathcal{O}_C(n_1,...,n_r))=i^*\mathcal{O}_{\mathbb{P}^n}(1).$
\end{proof}
Next, we recall that a $k$-{\it pointed chain} (of rational curves) is a datum
$(C,B_1,...,B_k)$ consisting of: (i) a chain $C=C_1\cup...\cup C_r$ of rational curves,
(ii) a set of $k$ distinct points $B_1,...,B_k\in C$ which for $r>1$ are also distinct
from the points $A_i=C_i\cap C_{i+1},\ i=1,...,r-1$. For any $k$-pointed chain
$(C,B_1,...,B_k)$ we denote its isomorphism class by $[(C,B_1,...,B_k)]$. The set
$\overline{M_{0,k}}$ of isomorphism classes of $k$-pointed chains is the well known
{moduli space of $k$-pointed chains (of rational curves)}.

Denote by $\Hilb^{n+1}\mathbb{P}^n$ the Hilbert scheme of subschemes of $\mathbb{P}^n$ with
Hilbert polynomial $n+1$. Fix $n+2$ points $B_1,...,B_{n+2}\in\mathbb{P}^n$ in general
position, i.e. such that no $n+1$ points lie in a hyperplane of $\mathbb{P}^n$.
Consider the set
$V(B_1,...,B_{n+2}):= \{C\in\Hilb^{n+1}\mathbb{P}^n| (C,B_1,...,B_k)$\ is\ a\ $k$-pointed chain
in\ $\mathbb{P}^n\}$ and the morphism
$$
\theta:V(B_1,...,B_{n+2})\to\overline{M_{0,n+2}}, C\mapsto[(C,B_1,...,B_{n+2})].
$$
In addition, put $\Gamma(B_1,...,B_{n+2}):=\{(x,C)\in\mathbb{P}^n\times
V(B_1,...,B_{n+2})|\ x\in C\}.$

Now we invoke results of Kapranov \cite{K} concerning $k$-pointed chains.
Together with \refle{embed} these results yield the following proposition.

\begin{proposition}\label{Kapranov}~

1) The morphism $\theta$ is an isomorphism, hence it induces an embedding
$i_\Gamma:\Gamma(B_1,...,B_{n+2})\hookrightarrow\mathbb{P}^n\times\overline{M_{0,n+2}}$.

2) For any strict partition $(n_1,...,n_r)$ of $n\in\ZZ_{>0}$ and any polarized chain
$(C_0,\mathcal{O}_{C_0}(n_1,...,n_r))$, there exist points $B_1,...,B_{n+2}\in C_0$ such
that the point $c_0=[(C_0,B_1,...,B_{n+2})]\in\overline{M_{0,n+2}}$ satisfies the condition
$\mathcal{O}_{C_0}(n_1,...,n_r)\simeq i_\Gamma^*(\mathcal{O}_{\mathbb{P}^n}(1)
\boxtimes\mathcal{O}_{\overline{M_{0,n+2}}})_{|C_0\times\{c_0\}})$.

3) The family of curves
$\pi_\Gamma:\Gamma(B_1,...,B_{n+2})\overset{i_\Gamma}\hookrightarrow\mathbb{P}^n\times
\overline{M_{0,n+2}}\overset{pr_2}\to\overline{M_{0,n+2}}$ is a semiuniversal deformation
of the curve $C_0$, i.e., for any flat family $\pi:X\to B$ of chains of rational curves
such that $C_0=\pi^{-1}(0)$ for some point $0\in B$, there exists a neighbourhood $U\ni0$
in $B$ and a morphism $f:U\to\overline{M_{0,n+2}}$ with $f(0)=c_0$ and
$\pi^{-1}(U)=\Gamma(B_1,...,B_{n+2})\times_{\overline{M_{0,n+2}}}U.$

4) The line bundle
$\mathcal{L}:=\Phi^*i_\Gamma^*(\mathcal{O}_{\mathbb{P}^n}(1)\boxtimes
\mathcal{O}_{\overline{M_{0,n+2}}})$,
where
$\Phi:\pi^{-1}(U)\to\Gamma(B_1,...,B_{n+2})$
is the induced morphism, satisfies the property
$\mathcal{L}_{|C_0}\simeq \mathcal{O}_{C_0}(n_1,...,n_r).$
\end{proposition}

\begin{corollary}\label{exists}
Let $\pi:X\to B$ be a flat family of chains of rational curves.
Let
$0\in B$ be a fixed point.
For any line bundle
$L_0$ on the fiber $C_0=\pi^{-1}(0)$
there exists
a neighbourhood $U\ni0$ in $B$ and a line bundle
$\mathcal{L}$
on $\pi^{-1}(U)$ such that
$\mathcal{L}_{|C_0}\simeq L_0$.
\end{corollary}
\begin{proof}
If $L_0$ is ample, its restrictions to each irreducible component of the fiber $C_0$ define a strict
partition $(n_1,...,n_r)$, and our statement is an immediate consequence of Proposition
\ref{Kapranov},2). Since any line bundle $L_0$ on $C_0$ can be represented as
$L_0'\otimes L_0''^\vee$ for some ample bundles $L_0',\ L_0''$, the Corollary follows.
\end{proof}

\begin{lemma}\label{ineq2}
Let $C=C_1\cup...\cup C_r$, be a chain of rational curves and let $E$ be a vector bundle
on $C$. Then, for any $\delta_1,...,\delta_r\in\ZZ_{\geq 0}$ with $\overset{r}{\underset{i=1}\sum}\delta_i>0,$ one has
\begin{equation*}
h^0(E(-\mathbf{d}_1(E_{|C_1})-\delta_1,...,-\mathbf{d}_1(E_{|C_r})-\delta_r))=0.
\end{equation*}
\end{lemma}
\begin{proof}
We use induction on $r$. For $r=1$ the statement is clear from the definition of
$\mathbf{d}_1(E)$. For the step of induction we just consider the case $r=2$, since for
arbitrary  $r$ the argument goes through without changes. Let $r=2$ and $\delta_1>0$. Then
clearly
$h^0(E(-\mathbf{d}_1(E_{|C_1})-\delta_1))=0$. Hence the natural exact triple
\begin{equation*}
0\to (E_{|C_2})(-\mathbf{d}_1(E_{|C_2})-\delta_2-1)\to
E(-\mathbf{d}_1(E_{|C_1})-\delta_1,-\mathbf{d}_1(E_{|C_2})-\delta_2)\to
\end{equation*}
$$
\to(E_{|C_1})(-\mathbf{d}_1(E_{|C_1})-\delta_1)\to0
$$
implies the equality
$h^0(E(-\mathbf{d}_1(E_{|C_1})-\delta_1,-\mathbf{d}_1(E_{|C_2})-\delta_2))=0.$
\end{proof}

\vspace{1.5cm}
\section{Construction of special rational curves in $G(r,V)$}
\label{rational curves}
\vspace{0.5cm}

Let $E$ be a rank $k$ vector bundle on $G=G(r,V)$. Denote
\begin{equation*}\label{D(E)}
\mathbf{D}(E):=\max\{\mathbf{D}(E_{|l})|\ l\ {\rm is\ a\ line\ in}\ G\}.
\end{equation*}
Our aim in this section is to prove that, for any
point $y_0\in G$, the inequality $\mathbf{D}(E_{|C})\le2r\mathbf{D}(E)$ holds on a dense
open subset of a suitably defined subscheme of the Hilbert scheme $H_{2r}(y_0)$ of rational
curves $C$ of degree $2r$ on $G$ passing through the point $y_0$.

We start with the following construction. Under the assumption that $\dim V\ge3r$,
let $V',V'',V'''$ be three $r$-dimensional subspaces of $V$ such that
\begin{equation}\label{V'}
V'\cap V''=V'\cap V'''=V''\cap V'''=\{0\},\ \ \ V'''\subset V'\oplus V''.
\end{equation}
In addition, fix $r$ linearly independent one-dimensional subspaces
$V_i',\ i=1,...,r,$ in $V'$. This datum defines linearly independent one-dimensional subspaces
$V_i'':=V''\cap(V'''\oplus V_i')$  in $V''$, as well as linearly independent one-
dimensional subspaces $V_i''':=V'''\cap(V_i'\oplus V_i'')$.
We obtain $r$ projective lines
$\mathbb{P}^1_i:=P(V_i'\oplus V_i'')$,
with points $V'_i,V''_i,V'''_i,\ i=1,...,r,$ on them.
On each of the lines $\mathbb{P}^1_i$ there is an affine coordinate $t_i$ uniquely determined by
the condition
\begin{equation*}\label{coord-t}
V'_i=\{t_i=0\},\ \ V''_i=\{t_i=\infty\},\ \ V'''_i=\{t_i=1\},\ \ \ \ i=1,...,r.
\end{equation*}
Let $V_{ti}\in\mathbb{P}^1_i$ be the point with affine coordinate $t$.
By construction, the points $V_{t1},...,V_{tr}$, considered as one-dimensional subspaces of $V$,
are linearly independent in $V$ and their
span $V^r_t:=V_{t1}\oplus...\oplus V_{tr}$
is an $r$-dimensional subspace of $V$. Thus we have an embedding
\begin{equation*}\label{embed-phi}
\phi_1:\ \mathbb{P}^1\hookrightarrow G,\ t\mapsto V^r_t
\end{equation*}
such that
$\phi_1(0)=V',\ \phi_1(\infty)=V'',\
\phi_1(1)=V''',\ \ \ \ \phi_1^*\mathcal{O}_G(1)=\mathcal{O}_{\mathbb{P}^1}(r)$.
Note that the degree $r$ curve $C^r_1={\rm im}\phi_1$ depends only on the choice of the
triple of $r$-dimensional spaces $V',V'',V'''$. We call the curve $C^r_1$ the {\it Segre curve}
associated to $V',V'',V'''$.
Moreover, the subspaces $V',V'',V'''$ define an embedding
\begin{equation*}\label{Segre embed}
s=s(V',V'',V'''):\mathbb{P}^1\times\mathbb{P}^{r-1}\overset{s_{1,1}}\hookrightarrow
\mathbb{P}^{2r-1}\overset{j}\hookrightarrow P(V),
\end{equation*}
where $s_{1,1}$ is the Segre embedding and $j$ is an embedding with ${\rm im}j=P(V'\oplus V'')$. By construction,
$\phi_1(t)=s(V',V'',V''')(\{t\}\times\mathbb{P}^{r-1}).$
We call $s(V',V'',V''')$ the {\it extended Segre embedding} associated to $V',V'',V'''$.

More generally, for any $(t_2,...,t_r)\in (\mathbf{k}^*)^{r-1}$ the triple of
spaces $V',V'',V^{(t_2,...,t_r)}$, where
\begin{equation*}\label{V^...}
V^{(t_2,...,t_r)}:=V_{11}\oplus V_{t_22}\oplus...\oplus V_{t_rr},
\end{equation*}
satisfies condition (\ref{V'}) with $V'''$ replaced by $V^{(t_2,...,t_r)}$
and hence yields a Segre curve
\begin{equation}\label{C^r t_i}
C^r_{(t_2,...,t_r)}:={\rm im}~s(V',V'',V^{(t_2,...,t_r)})
\end{equation}
in $G$. In particular, $C^r_{(1,...,1)}$ coincides with the Segre curve $C^r_1={\rm im}\phi_1$. In addition,
\begin{equation*}\label{C^r}
y_0:=\{V'\}\in C^r_{(t_2,...,t_r)}.
\end{equation*}

Let $l_1$ be the line in $G$ determined by the flag
$(V'_2\oplus..\oplus V'_r\subset V''_1\oplus V'_1\oplus..\oplus V'_r)$
and let $l_i$, for $i=2,...,r$, be the line in $G$ determined by the flag
$(V''_1\oplus...\oplus V''_{i-1}\oplus V'_{i+1}\oplus...\oplus V'_{r}\subset
V''_1\oplus...\oplus V''_{i}\oplus V'_i\oplus...\oplus V'_r)$.
These lines constitute a chain of rational curves
\begin{equation*}\label{r-chain}
C_0^r=l_1\cup...\cup l_r
\end{equation*}
in $G$.
Moreover, setting $(t_2,t_3,...,t_r)=(t,t^2,...,t^{r-1})$,
one easily proves the following lemma.
\begin{lemma}\label{deform to r-chain}
Let $G=G(r,V)$, with $\dim V\ge3r$. Consider the surface
$\mathcal{S}^\times:=
\{(x,t)\in G\times \mathbf{k}^\times|\ x\in C^r_t:=C^r_{(t,t^2...,t^{r-1})}\}$
with projection $\pi^\times:\ \mathcal{S}^\times\to \mathbf{k}^\times,\ (x,t)\mapsto t$.
Then the following statements hold.

(i) The fibers $C^r_t=(\pi^\times)^{-1}(t)$ are rational curves passing through the point
$y_0\in G$.

(ii) Let $\mathbf{k}^\times\hookrightarrow \mathbf{k}=\mathbb{A}^1$ be the standard inclusion
and let $\mathcal{S}$ be the closure of $\mathcal{S}^\times$ in $G\times \mathbf{k}$. The
extended projection $\pi:\ \mathcal{S}\to\mathbb{A}^1$ is a flat morphism and the natural
morphism $\nu:\ \mathcal{S}\to G$, $(x,t)\mapsto x$ is birational and gives an isomorphism
$\nu:\pi^{-1}(0)\overset{\sim}\to C_0^r$.
\end{lemma}

We note next that the condition $\dim V\ge3r$ in Lemma \ref{deform to r-chain} can be removed.
In fact, let $\dim V<3r$. By our assumption (see \refsec{sec2}), $r+2\le\dim V$. Fix a $3r$-dimensional vector space $\tilde{V}$ and let $V',V'',V'''$ be three $r$-dimensional
subspaces of $\tilde{V}$ satisfying the conditions (\ref{V'}). By performing the above
constructions for this datum we obtain a surface $\tilde{\mathcal{S}}$ with projections
$\mathbb{A}^1\overset{\tilde{\pi}}\leftarrow \tilde{\mathcal{S}}
\overset{\tilde{\nu}}\to\tilde{G}:=G(r,\tilde{V})$ as in Lemma \ref{deform to r-chain}. For
$y\in\tilde{\mathcal{S}}$, let $V_y$ be the $r$-dimensional subspace of $\tilde{V}$
corresponding to the point $\tilde{\nu}(y)$ in $\tilde{G}$. Since $\dim \mathcal{S}=2$, it
follows from the inequality $r+2\le\dim V$ that there exists a subspace $L$ of dimension
$3r-\dim V$ in
$\tilde{V}$ such that
\begin{equation}\label{empty}
L\cap(\underset{y\in\tilde{\mathcal{S}}}{\cup} V_y)=\{0\}.
\end{equation}
Fix an isomorphism $h:\tilde{V}/L\overset{\sim}\to V$ and consider the rational morphism
$\tilde{f}:\tilde{G}\dasharrow G$, \mbox{$V^r\mapsto h((V^r+L)\mod L))$}.
Then (\ref{empty}) implies that the morphism
$f:\tilde{\mathcal{S}}\to G\times\mathbb{A}^1,(x,t)\mapsto(\tilde{f}(x),t)$
is an embedding, i.e. that there exists an isomorphism
$g:\mathcal{S}:=f(\tilde{\mathcal{S}})\overset{\sim}\to\tilde{\SS}$ such that
$f\circ g=\id_\mathcal{S}$. Hence the surface $\mathcal{S}$ with its projections
$\pi:=\tilde{\pi}\circ g:\mathcal{S}\to\mathbb{A}^1$
and
$\nu:=f\circ\tilde{\nu}\circ g:\mathcal{S}\to G$
satisfies the assertion of Lemma \ref{deform to r-chain}.

Theorem \ref{theorem 5.1}, Lemma \ref{deform to r-chain} and this latest argument directly imply the following corollary.
\begin{corollary}\label{bound A}
Let $E$ be a rank $k$ vector bundle on the Grassmannian $G$. There exists an open subset
$U(E)$ of $\mathbf{k}^*$
such that, in the notation of Lemma \ref{deform to r-chain}, the inequality
\begin{equation}\label{key ineq}
\mathbf{D}(E_{|C^r_t})\le r\mathbf{D}(E)
\end{equation}
holds for any $t\in U(E)$.
\end{corollary}

Fix $t\in U(E).$ According to (\ref{C^r t_i}), the extended Segre embedding
$s_t:=s(V',V'',V^{(t,t^2,...,t^{r-1})}):
\mathbb{P}^1\times\mathbb{P}^{r-1}\hookrightarrow P(V)
$
induces an embedding
$\psi_t:\mathbb{P}^1\to G,\ u\mapsto s_t(\{u\}\times\mathbb{P}^{r-1})$,
such that
\begin{equation}\label{im psi_t}
{\rm im}\psi_t=\pi^{-1}(t)=C^r_t
\end{equation}
is the Segre curve from Lemma \ref{deform to r-chain}.
We will now construct another Segre curve ${C'}^r_t$ in $G$ such that $C^r_t\cup{C'}^r_t$
is a chain of rational curves (see (\ref{2-chain}) below).

For this, assume again temporarily that $\dim V\ge3r$. Set
$W'':=V'',\ W''_i:=V''_i,\ i=1,...,r,$ and choose two $r$-dimensional subspaces
$W',W'''$ in $V$ satisfying the conditions similar to (\ref{V'})
\begin{equation}\label{W'}
W'\cap W''=W'\cap W'''=W''\cap W'''=\{0\},\ \ \ W'''\subset W'\oplus W'',
\end{equation}
and the condition
\begin{equation}\label{cap=0}
W'\cap (V'\oplus V'')=\{0\}.
\end{equation}
We repeat the above construction for the datum $(W',W'',W''')$ instead of $(V',V'',V''')$.
First, there are uniquely defined linearly independent one-dimensional subspaces
$W_i':=W'\cap(W'''\oplus W''_i)$  in $W'$, as well as linearly independent one-
dimensional subspaces
$W_i''':=W'''\cap(W'\oplus W''_i)$  in $W'''$.
Furthermore, for $i=1,...,r$ there is a uniquely defined affine coordinate $t_i$
on the projective line $P(W'_i\oplus W''_i)$ such that
$W'_i=\{t=0\},\ \ W''_i=\{t=\infty\},\ \ W'''_i=\{t=1\}$. Denote by $W_{ti}$ the point on
$P(W'_i\oplus W''_i)$ with coordinate $t_i$. Finally, set
$W^{(t_2,...,t_r)}:=W_{11}\oplus W_{t_22}\oplus...\oplus W_{t_rr},
\ (t_2,...,t_r)\in (\mathbf{k}^*)^{r-1}$.

Next, for any $t\in\mathbf{k}^*$ consider the extended Segre embedding
$s'_t:=s(W',W'',W^{(t,t^2,...,t^{r-1})}):
{\mathbb{P}^1}'\times\mathbb{P}^{r-1}\hookrightarrow P(V)$,
where ${\mathbb{P}^1}'$ is a copy of $\mathbb{P}^1$.
This yields a map
$\psi'_t:{\mathbb{P}^1}'\to G,\ u\mapsto s'_t(\{u\}\times\mathbb{P}^{r-1})$,
such that
\begin{equation}\label{im psi'_t}
{C'}^r_t={\rm im}\psi'_t
\end{equation}
is a degree-$r$ Segre curve and, as in Corollary \ref{bound A}, there exists an open subset
$U'(E)$ of $\mathbf{k}^*$ such that
\begin{equation}\label{key ineq 2}
\mathbf{D}(E_{|{C'}^r_t})\le r\mathbf{D}(E),\ \ \ \ t\in U'(E).
\end{equation}
Moreover, for $t\in U_0(E):=U(E)\cap U'(E)$ we have $C^r_t\cap {C'}^r_t=\{V''\}$, i.e.
\begin{equation}\label{2-chain}
C^r_t\cup {C'}^r_t,\ \ \ t\in U_0(E),
\end{equation}
is a chain of rational curves in $G$ passing through the point $y_0=\{V'\}$.

Note that the assumption $\dim V\ge3r$ can again be removed. Indeed, consider a $3r$-dimensional space $\tilde{V}$ containing the
$r$-dimensional subspaces $V',V''=W'',V''',W',W'''$ satisfying (\ref{V'}), (\ref{W'}) and
(\ref{cap=0}). Set $\tilde{G}:=G(r,\tilde{V})$ and let $\tilde{\mathcal{S}}'$ be a surface in
$\tilde{G}\times\mathbb{A}^1$ defined by $W',W'',W'''$ in the same way as $\tilde{\mathcal{S}}$ was
defined by $V',V'',V'''$. Consider the subspace $L\subset\tilde V$ of dimension
$3r-\dim V$ satisfying the following condition similar to (\ref{empty}):
\begin{equation*}\label{empty2}
L\cap(\underset{y\in\tilde{\mathcal{S}}\cup\tilde{\mathcal{S}}'}{\cup} V_y)=\{0\}.
\end{equation*}
Then the argument following \refle{deform to r-chain} goes through without change, in
particular respective morphisms
$f':\tilde{\mathcal{S}}\cup\tilde{\mathcal{S}}'\hookrightarrow G\times\mathbb{A}^1$,
$\pi':\mathcal{S}\cup\mathcal{S}':={\rm im}f\to\mathbb{A}^1$ and
$\nu':\mathcal{S}\cup\mathcal{S}'\to G$ are defined and the chain (\ref{2-chain}) is the image in $G$ of the chain $(\pi')^{-1}(0)$ under $\nu'$.

As a next step we construct rational curves $C^{2r}$ of degree $2r$ in $G$ by deforming the chain (\ref{2-chain}).
\begin{proposition}\label{upper bound}
Let $E$ be a rank $k$ vector bundle on $G$ and $y_0\in G$ be an arbitrary point. There exists a
rational degree $2r$ curve $C$ in $G$ passing through the point $y_0$ such that:

(i)\ \ \ $\mathbf{D}(E_{|C})\le2r\mathbf{D}(E),$

(ii)\ \ \ $Q_{|C}\simeq (\mathcal{O}_{C}(2))^r,$ where $Q$ is the antitautological
bundle, i.e. the bundle dual to the tautological bundle on $G$.

\end{proposition}
\begin{proof} Fix $t\in U_0(E)$. Assume first that $\dim V\ge3r$.
Since $V''=W''$, the embeddings $s_t$ and $s'_t$ define an embedding
$Y\times \PP^{r-1}\to P(V)$ where $Y=\mathbb{P}^1\cup{\mathbb{P}^1}'$ as a reducible conic with
singular point $w_0=\mathbb{P}^1\cap{\mathbb{P}^1}'$ and a marked point
$z_0\in\mathbb{P}^1\subset Y$ such that
$\tilde{s}_t(\{w_0\}\times\mathbb{P}^{r-1})=P(V''),\ \
\tilde{s}_t(\{z_0\}\times\mathbb{P}^{r-1})=P(V')$. Note that, as a consequence of (\ref{cap=0}),
there exists an embedding of $g:Y\hookrightarrow\mathbb{P}^2$ such that $\tilde{s}_t$ fits in
the composition of maps
\begin{equation}\label{composn}
\tilde{s}_t:\ Y\times\mathbb{P}^{r-1}\overset{g\times id}\hookrightarrow
\mathbb{P}^2\times\mathbb{P}^{r-1}\overset{s_{1,1}}\hookrightarrow\mathbb{P}^{3r-1}
\overset{j}\hookrightarrow P(V),
\end{equation}
where $s_{1,1}$ is the Segre embedding by the linear series
$|\mathcal{O}_{\mathbb{P}^2}(1)\boxtimes\mathcal{O}_{\mathbb{P}^{r-1}}(1)|$
and $j$ is an embedding.

Now consider a pencil of conics
$\{Y_\tau\}_{\tau\in\mathbb{P}^1}$ in $\mathbb{P}^2$ satisfying the conditions:
(i) $Y_0=Y$,
(ii) all conics of the pencil pass through the point $z_0$ and
(iii) the generic conic in the pencil is smooth (i. e. a rational curve). Set
$\mathcal{U}'=\{\tau\in\mathbb{P}^1\smallsetminus\{0\} \ |\ Y_\tau$ is smooth$\}$ (this is a
dense open subset of $\mathbb{P}^1$). In view of (\ref{composn})
there exists a dense open subset $\mathcal{U}^*$ of $\mathcal{U}'$ such that, for any
$\tau\in\mathcal{U}^*$ the composition
$\tilde{s}_t:\ \mathbb{P}^1\times\mathbb{P}^{r-1}\simeq
Y_\tau\times\mathbb{P}^{r-1}\hookrightarrow
\mathbb{P}^2\times\mathbb{P}^{r-1}\overset{s_{1,1}}\hookrightarrow
\mathbb{P}^{3r-1}\overset{j}\hookrightarrow P(V),$
coincides with the embedding
$f_{t,\tau}:\ \mathbb{P}^1\times\mathbb{P}^{r-1}\hookrightarrow P(V)$
by a subseries of the linear series
$|\mathcal{O}_{\mathbb{P}^1}(2)\boxtimes\mathcal{O}_{\mathbb{P}^{r-1}}(1)|$.
This implies that the induced map
\begin{equation}\label{phi_t,tau}
\phi_{t,\tau}:\ \mathbb{P}^1\to G, \ u\mapsto f_{t,\tau}(\{u\}\times\mathbb{P}^{r-1})
\end{equation}
satisfies the property
\begin{equation}\label{rO(2)}
\phi_{t,\tau}^*Q\simeq (\mathcal{O}_{\mathbb{P}^1}(2))^r,\ \ \ \tau\in\mathcal{U}^*.
\end{equation}
Put $\mathcal{U}:=\mathcal{U}^*\cup\{0\}$
and consider the total space
$\Pi_t:=\{(x,\tau)\in\mathbb{P}^2\times\mathcal{U}|\ x\in Y_\tau\}$
of the above pencil of conics, together with the projections
$\mathbb{P}^2\overset{\sigma}\leftarrow\Pi_t\overset{\rho}\rightarrow\mathcal{U}.$
We obtain a morphism
\begin{equation*}\label{phi_t}
\phi_t:\ \Pi_t\to G, (x,\tau)\mapsto s_{1,1}(\{x\}\times\mathbb{P}^{r-1}).
\end{equation*}
By construction, $\rho^{-1}(\tau)=Y_\tau$ and for $\tau\in\mathcal{U}^*$
the map $\phi_{t|Y_\tau}$ coincides with $\phi_{t,\tau}$ from (\ref{phi_t,tau}). Moreover, by (\ref{im psi_t}) and (\ref{im psi'_t}) we have an isomorphism
\begin{equation}\label{phi_t|Y_0}
\phi_{t|Y_0}:\ Y_0\overset{\sim}\to C^{2r}_{t,0}:=C^r_t\cup{C'}^r_t.
\end{equation}
This means that $\phi_{t|Y_\tau}$ is an embedding near $\tau=0,$ i.e. the set
$\mathcal{U}_0=\{\tau\in\mathcal{U}|\ Y_\tau$ is smooth and $\phi_{t|Y_\tau}$ is an
embedding$\}$ is dense in $\mathcal{U}$. We thus obtain isomorphisms
\begin{equation}\label{phi_t,tau2}
\phi_{t,\tau}:\ \mathbb{P}^1\overset{\sim}\to C^{2r}_{t,\tau}={\rm im}\phi_{t,\tau}\subset G,\ \ \
\tau\in\mathcal{U}_0.
\end{equation}
The isomorphisms (\ref{phi_t|Y_0}) and (\ref{phi_t,tau2}) show that
$\{C^{2r}_{t,\tau}\}_{\tau\in\mathcal{U}_0\cup\{0\}}$
is a flat family of curves in $G$ whose fiber at 0 is a chain of rational curves of the form
$C^{2r}_{t,0}=C^r_t\cup{C'}^r_t$
and whose other fibers are rational curves $C^{2r}_{t,\tau}$. Hence, applying Theorem
\ref{theorem 5.1} to $\phi_{t,r}$ we obtain that
$\mathcal{U}(E):=\{\tau\in\mathcal{U}_0\cup\{0\}|\
\mathbf{D}(E_{|C^{2r}_{t,\tau}})\le\mathbf{D}(E_{|C^r_t})+\mathbf{D}(E_{|{C'}^r_t})\}$
is a dense open subset of $\mathcal{U}_0\cup\{0\}$.
Combining this with (\ref{key ineq}) and (\ref{key ineq 2}), and using (\ref{rO(2)}), we obtain
the assertion of the Proposition for any curve
$C:=C^{2r}_{t,\tau},\ (t,\tau)\in U_0(E)\times\mathcal{U}_0(E)$.

Finally, it remains to remove the assumption $\dim V\ge3r$. Let $\dim V<3r$. Take a space
$\tilde{V}$ of dimension $3r$ and choose its subspace $L$
of dimension $3r-\dim V$ satisfying the condition
\begin{equation*}\label{empty3}
P(L)\cap s_{1,1}(\mathbb{P}^2\times\mathbb{P}^{r-1})=\emptyset,
\end{equation*}
where $s_{1,1}$ is the Segre embedding defined in (\ref{composn}) and where the intersection
is taken in the space $\mathbb{P}^{3r-1}$ which is identified with $P(\tilde{V})$ in view of the
condition (\ref{cap=0}). (Note that $L$ always exists as $\dim V\geq r+2$.)
The rest of the argument goes through as in the remark preceeding the Proposition.
\end{proof}

We are now ready to discuss Hilbert schemes. Recall that any rational curve of given degree $k$ in $G$  can be considered as a point in the Hilbert scheme $\Hilb^{kt+1}G$. Set
\begin{equation}\label{H_k}
H_k:=\{C\in\Hilb^{kt+1}G\ |\ C\ {\rm is\ a\ rational\ curve\ of\
degree~}k\mathrm{~in}\ G\},
\end{equation}
$$
R_k=\{\phi:\mathbb{P}^1\to G\ |\ \phi\ {\rm is\ an\ embedding}\}.
$$
It is well known (see, e.g., \cite[Theorem 2.1]{St}) that $H_k$ is a smooth irreducible open
subset of $\Hilb^{kt+1}G$ and that the natural morphism
\begin{equation}\label{g_k}
g_k:R_k\to H_k, \ \phi\mapsto{\rm im}\phi
\end{equation}
is a principal $PGL(2)$-bundle. Next, consider the vector space
$\Hom(V^\vee,\mathbf{k}^r\otimes H^0(\mathcal{O}_{\mathbb{P}^1}(2)))$
and its dense open subset
\begin{equation}\label{W_S^surj}
W:=\{e\in \Hom(V^\vee,\mathbf{k}^r\otimes H^0(\mathcal{O}_{\mathbb{P}^1}(2)))
\ |\ {\rm the\ composition}\
\end{equation}
$$
\tilde{e}:
V^\vee\otimes\mathcal{O}_{\mathbb{P}^1}\overset{e\otimes id}\to
\mathbf{k}^r\otimes H^0(\mathcal{O}_{\mathbb{P}^1}(2))\otimes\mathcal{O}_{\mathbb{P}^1}
\overset{ev}\to\mathbf{k}^r\otimes \mathcal{O}_{\mathbb{P}^1}(2)\
{\rm is\ an~epimorphism}\},
$$
where $ev$ is the evaluation map. Let
$
\gamma:\ V^\vee\otimes\mathcal{O}_G\to Q
$
be the natural epimorphism.
By the universality of the Grassmannian $G$ any element
$e\in W$ defines a pair
\begin{equation}\label{phi_e}
(\ \phi_e:\mathbb{P}^1\to G,\ \ \chi_e:\phi_e^*Q\overset{\sim}\to
\mathbf{k}^r\otimes \mathcal{O}_{\mathbb{P}^1}(2)\ )
\end{equation}
such that
\begin{equation}\label{tilde e}
\chi_e\circ\phi_e^*\gamma=\tilde{e},
\end{equation}
where $\tilde{e}$ is defined in (\ref{W_S^surj}). Conversely, the element $e$ is recovered by
the pair $(\phi_e,\chi_e)$ since clearly $e$ is obtained from $\tilde{e}$ by passing to
sections:
\begin{equation}\label{h0(tilde e)}
e=H^0(\tilde{e}).
\end{equation}

Now put $k=2r$ in (\ref{H_k}) and consider the set
$$
H_{2r}^*:=\{C\in H_{2r}\ |\ \mathbf{D}_Q(C)=0,\ {\rm i.e.}\
Q_{|C}\simeq\mathbf{k}^r\otimes \mathcal{O}_{\mathbb{P}^1}(2) \}.
$$
By semicontinuity,
$H_{2r}^*$ is an open subset of $H_{2r}$.
Moreover, $H_{2r}^*$ is nonempty (and hence dense in $H_{2r}$)
since it contains
all curves
$C_{t,\tau}^{2r}$ from Proposition \ref{upper bound}.

\begin{theorem}\label{B_S(E,y_0)}
Fix a point $y_0\in G$ and put
\begin{equation}\label{H_2r(y0)}
H_{2r}^*(y_0):=\{C\in H_{2r}^*|y_0\in C\}.
\end{equation}
For any vector bundle $E$ on $G$, $B(E,y_0):=\{C\in H_{2r}^*(y_0)\ |\ \mathbf{D}(E_{C})\le2r\mathbf{D}(E)\}$
is a dense open subset of the irreducible variety $H_{2r}^*(y_0)$.
\end{theorem}
\begin{proof}
Since $H_{2r}^*$ is smooth and irreducible, and the group $PGL(V^\vee)$ acts transitively
on $G$, it follows that $H_{2r}^*(y_0)$ is an irreducible (and smooth) subvariety of $H_{2r}^*$
which contains the curve $C=C_{t,\tau}^{2r}$ from Proposition \ref{upper bound}.
Moreover, since the condition
$\mathbf{D}(E_{C})\le2r\mathbf{D}(E)$ is open on
$C\in H_{2r}^*(y_0)$ by semicontinuity, Proposition \ref{upper bound} immediately implies the Theorem.
\end{proof}

Now take an arbitrary curve $C\in H_{2r}^*$ and pick an embedding
$\phi_C:\mathbb{P}^1\hookrightarrow G$ such that ${\rm im}\phi_C=C$.
In addition, pick an isomorphism $\chi_C:\phi_C^*Q\overset{\sim}\to
\mathbf{k}^r\otimes \mathcal{O}_{\mathbb{P}^1}(2)$. These data define an element
$e=H^0(\chi_C\circ\phi_C^*\gamma)\in W$
(cf. (\ref{tilde e}) and (\ref{h0(tilde e)})) such that $\phi_e=\phi_C$. Moreover, $e$ belongs to the subset
\begin{equation}\label{W_S^*}
W^*:=\{e\in W\ |\ \phi_e:\mathbb{P}^1\to G\ {\rm is\ an\ embedding}\}
\end{equation}
of $W$. This nonempty subset is clearly open in $W$, hence it is dense in $W$. Moreover,
setting $R_{2r}^*:=g_{2r}^{-1}(H_{2r}^*)$ (see (\ref{g_k})), we obtain a principal $GL(r)$-bundle $\theta_{2r}:W^*\to R_{2r}^*,\ e\mapsto\phi_e$. Since (\ref{g_k}) is a principal $PGL(2)$-bundle, the composition
\begin{equation*}\label{w_2r}
w_{2r}:=g_{2r}\circ\theta_{2r}:W^*\to H_{2r}^*
\end{equation*}
is a smooth surjective morphism.

Fix a point
$y_0\in G$
and consider the set
\begin{equation}\label{W_S^*(E,y_0)}
W^*(E,y_0):=w_{2r}^{-1}(B(E,y_0)).
\end{equation}
By (\ref{phi_e}), (\ref{W_S^*}) and (\ref{H_2r(y0)}) any point $e\in W^*(E,y_0)$ defines an
embedding $\phi_e:\mathbb{P}^1\hookrightarrow G$ with $y_0\in{\rm im}\phi_e$.
Let $z:=\phi_e^{-1}(y_0)$. We have an exact triple
$0\to\mathbf{k}^r\otimes\mathcal{O}_{\mathbb{P}^1}(1)\to
\mathbf{k}^r\otimes\mathcal{O}_{\mathbb{P}^1}(2)\to\mathbf{k}^r\otimes\mathbf{k}(z)\to0$.
By applying $\mathcal{H}om(V^\vee\otimes\mathcal{O}_{\mathbb{P}^1},-)$
and passing to sections we obtain an exact triple
$$
0\to\Hom(V^\vee,\mathbf{k}^r\otimes H^0(\mathcal{O}_{\mathbb{P}^1}(1)))\to
\Hom(V^\vee,\mathbf{k}^r\otimes H^0(\mathcal{O}_{\mathbb{P}^1}(2)))
\overset{res_z}\to\Hom(V^\vee,\mathbf{k}^r\otimes\mathbf{k}(z))\to0.
$$
By construction, the set
\begin{equation*}\label{W_S(E,y_0)}
W(E,y_0):=res_z^{-1}(res_z(e))\simeq
\Hom(V^\vee,\mathbf{k}^r\otimes H^0(\mathcal{O}_{\mathbb{P}^1}(1)))
\end{equation*}
depends only on $E$ and $y_0$ and contains $W^*(E,y_0)$ as a dense open subset. We thus obtain
the following corollary.
\begin{corollary}\label{cor B_S}
$W^*(E,y_0)=w_{2r}^{-1}(B(E,y_0))$ is a dense open subset of
$W(E,y_0)\simeq\Hom(V^\vee,\mathbf{k}^r\otimes H^0(\mathcal{O}_{\mathbb{P}^1}(1)))$.
\end{corollary}

\vspace{0.5cm}

\section{Construction of a special vector bundle on $\PP^1\times\PP^1$}
\label{vector bdl}
\vspace{0.5cm}

A key ingredient in the proof of our main result, Theorem \ref{main}, is a specific morphism of
$\mathbb{P}^1\times\mathbb{P}^1$
into $G=G(r,V)$. It defines a vector bundle of rank $r$ on
$\mathbb{P}^1\times\mathbb{P}^1$ as the pullback of the antitautological bundle $Q$ on $G$.
We construct this morphism in two steps. We first construct a
suitable vector bundle on
$\mathbb{P}^1\times\mathbb{P}^1$ and then prove that this bundle arises from an
appropriate morphism
of $\mathbb{P}^1\times\mathbb{P}^1$ to $G$.

We now proceed to the construction of a special vector bundle  on
$X:=\mathbb{P}^1\times\mathbb{P}^1$.
Let $pr_i:X\rightarrow\mathbb{P}^1,\ i=1,2,$ be the natural projections, and let
$0$ and $\infty$ be two fixed points on $\mathbb{P}^1$. Set
$P:=pr_1^{-1}(0),\ \ P':=pr_1^{-1}(\infty)$.

Fix a partition $(a_1,a_2,...,a_r)$ of $d\in\mathbb{Z}_{\geq 0}$
such that $a_1\ge a_2\ge...\ge a_r$,
and consider a set $Z'$ of $d-a_1$ distinct points on
$\mathbb{P}^1$ of the form
\begin{equation*}\label{Z'}
Z'=\overset{r}{\underset{i=2}\sqcup}Z'_i,\ \ \ \ Z'_i=\overset{a_i}{\underset{j=1}\sqcup}x_{ij}
\end{equation*}
($Z_i'\neq \varnothing$ for $a_i=0$). The set $Z'$ determines the set
\begin{equation*}\label{Z}
Z:=P'\cap pr_2^{-1}(Z')=\overset{r}{\underset{i=2}\sqcup}Z_i,\ \ \ \
Z_i:=P'\cap pr_2^{-1}(Z'_i)=\overset{a_i}{\underset{j=1}\sqcup}(\infty,x_{ij}).
\end{equation*}
In what follows we think of $Z_i$ as reduced 0-dimensional subschemes of $X$.

Setting
\begin{equation}\label{E_1}
\mathcal{E}_1:=\mathcal{O}_X(0,a_1)
\end{equation}
we will recursively construct sheaves $\mathcal{E}_k$ for $2\le k\le r$
via the exact triples
\begin{equation}\label{E_k}
0\to\mathcal{E}_{k-1}\to\mathcal{E}_k\to\mathcal{I}_{Z_k,X}(0,a_k)\to0,\ \ \ 2\le k\le r.
\end{equation}

\begin{proposition}\label{quadric2}
Assume that for $2\le k\le r$ the sheaf $\mathcal{E}_{k-1}$ is locally
free of rank $k-1$ and satisfies the conditions
\begin{equation}\label{H^2=0}
h^2(\mathcal{E}_{k-1}(0,-a_i))=0,\ \ \ k\le i\le r,
\end{equation}
\begin{equation}\label{some h^2=0}
h^2(\mathcal{E}_{k-1}(0,-1))
=0,
\end{equation}
\begin{equation}\label{some h^1=0}
h^1(\mathcal{E}_{k-1})=h^1(\mathcal{E}_{k-1}(0,-1))=h^1(\mathcal{E}_{k-1}(1,0))=0,
\end{equation}
\begin{equation}\label{2 some h^1=0}
h^1(\mathcal{E}_{k-1}(1,-1))
=0,
\end{equation}
\begin{equation}\label{some h^0}
h^0(\mathcal{E}_{k-1}(0,-1))=a_1,\ \ h^0(\mathcal{E}_{k-1}(1,-1))=2a_1+a_2+...+a_{k-1},\ \
h^0(\mathcal{E}_{k-1})=a_1+k-1,\ \
\end{equation}
\begin{equation}\label{2 some h^0}
h^0(\mathcal{E}_{k-1}(1,0))=2a_1+a_2+...+a_{k-1}+2(k-1),\ \
h^0(\mathcal{E}_{k-1}(2,0))=3a_1+2a_2+...+2a_{k-1}+3(k-1),
\end{equation}
\begin{equation}\label{res E_k-1}
{\mathcal{E}_{k-1}}_{|pr_2^{-1}(x)}\simeq(\mathcal{O}_{\mathbb{P}^1})^{k-1},\ \
x\not\in Z'_2\sqcup...\sqcup Z'_{k-1},\ \ \ k\ge3,
\end{equation}
\begin{equation}\label{2 res E_k-1}
\mathcal{E}_{{k-1}|pr_2^{-1}(x)}\simeq(\mathcal{O}_{\mathbb{P}^1})^{k-3}
\oplus\mathcal{O}_{\mathbb{P}^1}(1)\oplus\mathcal{O}_{\mathbb{P}^1}(-1),\ \
x\in Z'_2\sqcup...\sqcup Z'_{k-1},\ \ \ k\ge3,
\end{equation}
\begin{equation}\label{E_k-1|P}
\mathcal{E}_{{k-1}|P}\simeq\mathcal{O}_{\mathbb{P}^1}(a_1)
\oplus...\oplus\mathcal{O}_{\mathbb{P}^1}(a_{k-1}),
\end{equation}
\begin{equation}\label{E_k-1|P'}
\mathcal{E}_{{k-1}|P'}\simeq\mathcal{O}_{\mathbb{P}^1}(a_1+...+a_{k-1})
\oplus(\mathcal{O}_{\mathbb{P}^1})^{k-2}.
\end{equation}
Then

(i) there is an epimorphism
\begin{equation}\label{Exts}
\Ext^1(\mathcal{I}_{Z_k,X}(0,a_k),\mathcal{E}_{k-1})\overset{\beta}\twoheadrightarrow
H^0(\EXT^1(\mathcal{I}_{Z_k,X}(0,a_k),\mathcal{E}_{k-1})),
\end{equation}
and moreover
\begin{equation}\label{eq64}
H^0(\EXT^1(\mathcal{I}_{Z_k,X}(0,a_k),\mathcal{E}_{k-1}))\simeq
H^0(\mathcal{O}_{Z_k}^{k-1})\simeq\underset{x\in Z_k}\oplus\mathbf{k}(x)^{k-1};
\end{equation}
(ii) there exists an element $\xi\in\Ext^1(\mathcal{I}_{Z_k,X}(0,a_k),\mathcal{E}_{k-1})$
such that the sheaf $\mathcal{E}_k$ defined by the corresponding exact triple (\ref{E_k})
is locally free of rank $k$ and satisfies the conditions
(\ref{H^2=0})-(\ref{E_k-1|P'}) with $k$ substituted for $k-1$; we label the so modified
conditions as (\ref{H^2=0}')-(\ref{E_k-1|P'}').
\end{proposition}

\begin{proof}
(i) The existence of the epimorphism (\ref{Exts}) follows from the standard exact sequence of
local and global Ext's
$$
0\to H^1(\mathcal{H}om(\mathcal{I}_{Z_k,X}(0,a_k),\mathcal{E}_{k-1}))\to
\Ext^1(\mathcal{I}_{Z_k,X}(0,a_k),\mathcal{E}_{k-1})\to
H^0(\EXT^1(\mathcal{I}_{Z_k,X}(0,a_k),\mathcal{E}_{k-1}))\to
$$
$$
\to
H^2(\mathcal{H}om(\mathcal{I}_{Z_k,X}(0,a_k),\mathcal{E}_{k-1}))
$$
and from (\ref{H^2=0}) in view of the canonical isomorphism
$\mathcal{H}om(\mathcal{I}_{Z_k,X}(0,a_k),\mathcal{E}_{k-1})=\mathcal{E}_{k-1}(0,-a_k)$.
The isomorphisms in (\ref{eq64}) are standard.

(ii) Pick an element $\xi\in\Ext^1(\mathcal{I}_{Z_k,X}(0,a_k),\mathcal{E}_{k-1})$ and consider
the extension (\ref{E_k}) defined by $\xi$.
Put $S_x:=pr_2^{-1}(x).$ If
$x\not\in Z'_2\sqcup...\sqcup Z'_{k-1},$ then
in view of (\ref{res E_k-1}) the restriction of (\ref{E_k}) onto $S_x$ is
$0\to(\mathcal{O}_{\mathbb{P}^1})^{k-1}\to\mathcal{E}_{k|S_x}\to\mathcal{O}_{\mathbb{P}^1}\to0$,
i.e.
$\mathcal{E}_{k|S_x}\simeq(\mathcal{O}_{\mathbb{P}^1})^{k}$. This implies (\ref{res E_k-1}$'$).

If
$x\in Z'_2\sqcup...\sqcup Z'_{k-1},$
then in view of (\ref{2 res E_k-1}) the restriction of (\ref{E_k}) onto $S_x$ is
$0\to(\mathcal{O}_{\mathbb{P}^1})^{k-3}
\oplus\mathcal{O}_{\mathbb{P}^1}(1)\oplus\mathcal{O}_{\mathbb{P}^1}(-1)
\to\mathcal{E}_{k|S_x}\to\mathcal{O}_{\mathbb{P}^1}\to0$, i.e.
$\mathcal{E}_{k|S_x}\simeq(\mathcal{O}_{\mathbb{P}^1})^{k-2}
\oplus\mathcal{O}_{\mathbb{P}^1}(1)\oplus\mathcal{O}_{\mathbb{P}^1}(-1),$ and we obtain (\ref{2 res E_k-1}$'$) for $x\in Z'_2\sqcup...\sqcup Z'_{k-1}$.

For $x\in Z'_k,$ one has $\mathcal{I}_{Z_k,X}(0,a_k)_{|S_x}\simeq\mathbf{k}(\bar x)\oplus\mathcal{O}_{\mathbb{P}^1}(-1),$
where $\bar x:=(\infty,x)$. Therefore (\ref{res E_k-1}) yields an exact sequence
\begin{equation}\label{E_k|S_x}
0\to(\mathcal{O}_{\mathbb{P}^1})^{k-1}\to\mathcal{E}_{k|S_x}\to
\mathbf{k}(\bar x)\oplus\mathcal{O}_{\mathbb{P}^1}(-1)\to 0.
\end{equation}
Here the extension (\ref{E_k|S_x}) is given by an element
\begin{equation}\label{xi^1_x}
\xi^1_x\in\Ext^1(\mathcal{I}_{Z_k,X}(0,a_k)_{|S_x},\mathcal{E}_{{k-1}|S_x})\simeq
\Ext^1(\mathbf{k}(\bar x)\oplus\mathcal{O}_{\mathbb{P}^1}(-1),(\mathcal{O}_{\mathbb{P}^1})^{k-1})\simeq
\end{equation}
$$
\simeq H^0(\EXT^1(\mathbf{k}(\bar x)\oplus\mathcal{O}_{\mathbb{P}^1}(-1),(\mathcal{O}_{\mathbb{P}^1})^{k-1}))\simeq\mathbf{k}(\bar x)^{k-1}.
$$
Note that a sufficient condition for (\ref{2 res E_k-1}$'$) is that $\xi_x^1\neq 0$ for $x\in Z'_k$. Note in addition that the restriction of (\ref{E_k}) onto $S_x$ defines a natural restriction map
$$
\psi_x:\Ext^1(\mathcal{I}_{Z_k,X}(0,a_k),\mathcal{E}_{k-1})\to
\Ext^1(\mathcal{I}_{Z_k,X}(0,a_k)_{|S_x},\mathcal{E}_{{k-1}|S_x})
$$
such that
$$
\xi^1_x=\psi_x(\xi).
$$
The map $\psi_x$ together with (\ref{Exts}) and (\ref{xi^1_x}) fits in the diagram
\begin{equation}\label{diag1}
\xymatrix{
\Ext^1(\mathcal{I}_{Z_k,X}(0,a_k),\mathcal{E}_{k-1})\ar@{->>}[d]_{\beta}
\ar[r]^{\psi_x\ \ \ }&
\Ext^1(\mathcal{I}_{Z_k,X}(0,a_k)_{|S_x},\mathcal{E}_{{k-1}|S_x})\ar[d]_{\simeq}\\
H^0(\mathcal{O}_{Z_k}^{k-1})\ar@{->>}[r]^{res_x}& \mathbf{k}(\bar x)^{k-1},
}
\end{equation}
where $res_x$ is the restriction epimorphism defined by the inclusion
$\bar x\hookrightarrow Z_k$.

Next, (\ref{E_k}) and (\ref{E_k-1|P}) give
$$
0\to\mathcal{O}_{\mathbb{P}^1}(a_1)\oplus...\oplus\mathcal{O}_{\mathbb{P}^1}(a_{k-1})
\to\mathcal{E}_{k|P}\to\mathcal{O}_{\mathbb{P}^1}(a_k)\to0
$$
Since $a_1\ge a_2\ge...\ge a_k,$ this extension splits and yields
(\ref{E_k-1|P}$'$). Furthermore, by (\ref{E_k}) and (\ref{E_k-1|P'}) we have
\begin{equation}\label{2 E_k|P'}
0\to\mathcal{O}_{\mathbb{P}^1}(a_1+...+a_{k-1})\oplus(\mathcal{O}_{\mathbb{P}^1})^{k-2}\to
\mathcal{E}_{k|P'}\to\mathcal{O}_{Z_k}\oplus\mathcal{O}_{\mathbb{P}^1}\to0.
\end{equation}
The extension (\ref{2 E_k|P'}) is given by an element
\begin{equation}\label{xi'}
\xi'\in\Ext^1(\mathcal{I}_{Z_k,X}(0,a_k)_{|P'},\mathcal{E}_{{k-1}|P'})\simeq
\Ext^1(\mathcal{O}_{Z_k}\oplus\mathcal{O}_{\mathbb{P}^1},
\mathcal{O}_{\mathbb{P}^1}(a_1+...+a_{k-1})\oplus(\mathcal{O}_{\mathbb{P}^1})^{k-2})\simeq
\end{equation}
$$
\simeq H^0(\mathcal{O}_{Z_k}^{k-1}).
$$
Since $a_1+...+a_{k-1}\ge0,$ there is a distinguished injection
$$
g_k=\underset{x\in Z_k}\oplus g_k(x):\ \underset{x\in Z_k}\oplus\mathbf{k}(x)\simeq
\Ext^1(\mathcal{O}_{Z_k}\oplus\mathcal{O}_{\mathbb{P}^1},
\mathcal{O}_{\mathbb{P}^1}(a_1+...+a_{k-1}))\hookrightarrow
\Ext^1(\mathcal{I}_{Z_k,X}(0,a_k)_{|P'},\mathcal{E}_{{k-1}|P'})
$$
$$
\simeq\underset{x\in Z_k}\oplus\mathbf{k}(x)^{k-1}.
$$
Furthermore, in view of (\ref{Exts}) and (\ref{xi'}), we have a diagram of morphisms similar to (\ref{diag1})
\begin{equation}\label{diag2}
\xymatrix{
\Ext^1(\mathcal{I}_{Z_k,X}(0,a_k),\mathcal{E}_{k-1})
\ar@{->>}[d]_{\beta}\ar@{->>}[r]^{\psi_{P'}\ \ \ }&
\Ext^1(\mathcal{I}_{Z_k,X}(0,a_k)_{|P'},\mathcal{E}_{{k-1}|P'})\ar[d]_{\simeq}\\
H^0(\mathcal{O}_{Z_k}^{k-1})\ar@{=}[r]& H^0(\mathcal{O}_{Z_k}^{k-1}).
}
\end{equation}
The diagrams (\ref{diag1}) and (\ref{diag2}) immediately imply that
the element
$\xi\in\Ext^1(\mathcal{I}_{Z_k,X}(0,a_k),\mathcal{E}_{k-1})$
can be chosen so that:

1) $(res_x\circ\beta)(\xi)\ne0$ for any $x\in Z_k$,

2) the element $\xi'$ in (\ref{xi'}) satisfies the condition
$\xi'=\psi_{P'}(\xi)\in{\rm im}g_k$.

\noindent It follows from these conditions that (\ref{2 res E_k-1}$'$) holds and the extension (\ref{2 E_k|P'}) implies
(\ref{E_k-1|P'}$'$).

To prove the remaining equalities (\ref{H^2=0}$'$)-(\ref{2 some h^0}$'$) for the vector bundle
$\mathcal{E}_k$ defined by $\xi$ as the extension (\ref{E_k}), we consider the standard
Koszul resolution
\begin{equation}\label{resoln}
0\to\mathcal{O}_X(-1,0)\to\mathcal{O}_X\oplus\mathcal{O}_X(-1,a_k)
\to\mathcal{I}_{Z_k,X}(0,a_k)\to0.
\end{equation}
Twisting (\ref{resoln}) by $\mathcal{O}_X(a,b)$ for appropriate $a,b$ and keeping in mind that $a_k\ge0$, we obtain
\begin{equation}\label{some h^2=0 I_Z}
h^2(\mathcal{I}_{Z_k,X}(0,a_k-a_i))=0,\ \ \ k+1\le i\le r,
\end{equation}
\begin{equation}\label{2 some h^2=0 I_Z}
h^2(\mathcal{I}_{Z_k,X}(0,a_k-1))
=0,
\end{equation}
\begin{equation}\label{some h^1=0 I_Z}
h^1(\mathcal{I}_{Z_k,X}(0,a_k))=
h^1(\mathcal{I}_{Z_k,X}(1,a_k))=0,
\end{equation}
\begin{equation}\label{2 some h^1=0 I_Z}
h^1(\mathcal{I}_{Z_k,X}(1,a_k-1))
=0,
\end{equation}
\begin{equation}\label{some h^0 I_Z}
h^0(\mathcal{I}_{Z_k,X}(0,a_k-1))=0,\ \
h^0(\mathcal{I}_{Z_k,X}(1,a_k-1))=a_k,\ \
h^0(\mathcal{I}_{Z_k,X}(0,a_k))=1,\ \
\end{equation}
\begin{equation}\label{2 some h^0 I_Z}
h^0(\mathcal{I}_{Z_k,X}(1,a_k))=a_k+2,\ \
h^0(\mathcal{I}_{Z_k,X}(1,a_k))=2a_k+3.
\end{equation}
Furthermore, twisting (\ref{E_k}) by $\mathcal{O}_X(a,b)$ we quickly see that:
(\ref{H^2=0}) and (\ref{some h^2=0 I_Z}) imply (\ref{H^2=0}$'$); (\ref{some h^2=0})
and (\ref{2 some h^2=0 I_Z}) imply (\ref{some h^2=0}$'$); (\ref{some h^1=0}) and
(\ref{some h^1=0 I_Z}) imply (\ref{some h^1=0}$'$); (\ref{2 some h^1=0}) and
(\ref{2 some h^1=0 I_Z}) imply (\ref{2 some h^1=0}$'$); (\ref{some h^0}) and
(\ref{some h^0 I_Z}) imply (\ref{some h^0}$'$); (\ref{2 some h^0}) and
(\ref{2 some h^0 I_Z}) imply (\ref{2 some h^0}$'$).
The Proposition is proved.
\end{proof}

As a corollary of Proposition \ref{quadric2} one obtains the following theorem.

\begin{theorem}\label{F}
For any $r\in\ZZ_{>0}$, $r\geq 2$, there exists a rank $r$ vector bundle $\mathcal{F}$ on
the surface
$X=\mathbb{P}^1\times\mathbb{P}^1$
with the following properties:

\begin{equation}\label{res F}
\mathcal{F}_{|pr_2^{-1}(x)}=(\mathcal{O}_{\mathbb{P}^1}(2))^r,\ \ x\not\in Z',
\end{equation}
\begin{equation}\label{2 res F}
\mathcal{F}_{|pr_2^{-1}(x)}=(\mathcal{O}_{\mathbb{P}^1}(2))^{r-2}
\oplus\mathcal{O}_{\mathbb{P}^1}(3)\oplus\mathcal{O}_{\mathbb{P}^1}(1),\ \ x\in Z',
\end{equation}
\begin{equation}\label{F|P}
\mathcal{F}_{|P}=\mathcal{O}_{\mathbb{P}^1}(a_1)\oplus...\oplus\mathcal{O}_{\mathbb{P}^1}(a_r),
\end{equation}
\begin{equation}\label{F|P'}
\mathcal{F}_{|P'}=\mathcal{O}_{\mathbb{P}^1}(d)\oplus(\mathcal{O}_{\mathbb{P}^1})^{r-1},
\end{equation}
\begin{equation}\label{h^1F}
h^1(\mathcal{F}(-1,-1))= h^1(\mathcal{F}(-2,0))=h^1(\mathcal{F}(-2,-1))=h^1(\mathcal{F}(-1,0))=0,
\end{equation}
\begin{equation}\label{h^2F}
h^2(\mathcal{F}(-2,-1))=0,
\end{equation}
\begin{equation}\label{3 some h^0}
h^0(\mathcal{F}(-1,-1))=a_1+d,\ \
h^0(\mathcal{F}(-1,0))=a_1+d+2r,\ \
h^0(\mathcal{F})=a_1+2d+3r.
\end{equation}
\end{theorem}

\begin{proof}
We define $\FF$ as $\mathcal{E}_r(2,0)$. The equalities (\ref{res F})-(\ref{3 some h^0}) follow directly from the equalities (\ref{E_1}) and (\ref{H^2=0}$'$)-(\ref{E_k-1|P'}$'$) for $k=r$ and from the observation that $a_1+...+a_r=d$.
\end{proof}

\section{Construction of a special morphism $f:\PP^1\times\PP^1\to G(r,V)$}
We are now ready to proceed with the construction of the desired morphism $f:X\to G$.
Fix a line $l_0$ in $G$, a point $y_0\in l_0$ and a degree $d$ morphism
$\psi':\mathbb{P}^1\to l_0$, and let
\begin{equation*}\label{psi}
\psi:\ \mathbb{P}^1\overset{\psi'}\to l_0\hookrightarrow G
\end{equation*}
be the composition.
The restriction of the canonical epimorphism
$\gamma:V^\vee\otimes\mathcal{O}_G\to Q$
to $l_0$ has the form
$g:V^\vee\otimes\mathcal{O}_{\mathbb{P}^1}\to
\mathcal{O}_{\mathbb{P}^1}(1)\oplus(\mathcal{O}_{\mathbb{P}^1})^{r-1}.$
Hence the epimorphism $\psi^*g$ has the form
\begin{equation*}\label{psi^*g}
\psi^*g:V^\vee\otimes\mathcal{O}_{\mathbb{P}^1}\twoheadrightarrow
\mathcal{O}_{\mathbb{P}^1}(d)\oplus(\mathcal{O}_{\mathbb{P}^1})^{r-1}.
\end{equation*}
Passing to sections
we obtain an element
\begin{equation}\label{g_psi}
g_\psi\in\Hom(V^\vee,
H^0(\mathcal{O}_{\mathbb{P}^1}(d)\oplus(\mathcal{O}_{\mathbb{P}^1})^{r-1})).
\end{equation}
Note that, similarly to (\ref{phi_e}), $\psi$ is determined by the element $g_\psi$.

Next, we put $k=d$ in (\ref{H_k}) and we fix a curve $C_0\in H_d$ together with an isomorphism
$\theta:\mathbb{P}^1\overset{\sim}\to C_0$.
The composition
$\phi_{C_0}:\mathbb{P}^1\overset{\theta}{\overset{\sim}\to} C_0\hookrightarrow G$
defines an epimorphism $\phi_{C_0}^*\gamma:V^\vee\otimes
\mathcal{O}_{\mathbb{P}^1}\twoheadrightarrow\phi_{C_0}^*Q$.
Moreover,
$$
\phi_{C_0}^*Q\simeq\mathcal{O}_{\mathbb{P}^1}(a_1)
\oplus...\oplus\mathcal{O}_{\mathbb{P}^1}(a_r)
$$
for some partition $(a_1,...,a_r)$ of $d$. (The nonnegativity of the integers $a_i$
follows from the surjectivity of
$\phi_{C_0}^*\gamma$.)
Pick an isomorphism
$\chi_{C_0}:\phi_{C_0}^*Q\overset{\sim}\to\mathcal{O}_{\mathbb{P}^1}(a_1)
\oplus...\oplus\mathcal{O}_{\mathbb{P}^1}(a_r)$.
The datum $(\phi_{C_0},\chi_{C_0})$ defines an element
\begin{equation*}\label{e(phi_C_0)}
e(C_0)=H^0(\chi_{C_0}\circ\phi_{C_0}^*\gamma)\in W(C_0)
:=\Hom(V^\vee,H^0(\mathcal{O}_{\mathbb{P}^1}(a_1)
\oplus...\oplus\mathcal{O}_{\mathbb{P}^1}(a_r)))
\end{equation*}
(cf. (\ref{tilde e}) and (\ref{h0(tilde e)})).

The set
$W(C_0)^{epi}:=\{e\in W(C_0)\ |\ {\rm the\ composition}\
V^\vee\otimes\mathcal{O}_{\mathbb{P}^1}\overset{e\otimes
id}\to
H^0(\mathcal{O}_{\mathbb{P}^1}(a_1)\oplus...\oplus\mathcal{O}_{\mathbb{P}^1}(a_r))
\otimes\mathcal{O}_{\mathbb{P}^1}
\overset{ev}\to\mathcal{O}_{\mathbb{P}^1}(a_1)\oplus...\oplus\mathcal{O}_{\mathbb{P}^1}(a_r)\
{\rm is\ an~epimorphism}\}$
is a dense open subset of $W(C_0)$ containing $e(C_0)$, and (by the universality property of $G$) any element $e\in
W(C_0)^{epi}$ determines a morphism
\begin{equation*}\label{phi_e2}
\phi_e:\mathbb{P}^1\to G
\end{equation*}
(cf. (\ref{phi_e})). For any vector bundle $E$ on $G$ we put
\begin{equation}\label{B_P(E,C_0)}
B_P(E,C_0):=\{C\in H_d\ |\ \mathbf{D}(E_{|C})\le\mathbf{D}(E_{|C_0})\}.
\end{equation}
By construction $C_0\in B_P(E,C_0)$, and moreover by semicontinuity, $B_P(E,C_0)$ is an open
subset of $H_d$. Since $e(C_0)\in W(C_0)^{epi}$ it follows that
\begin{equation}\label{W_P(E,C_0)}
W_P(E,C_0):=\{e\in W(C_0)^{epi}\ |\ {\rm im}\phi_e\in B_P(E,C_0)\}
\end{equation}
is a dense open subset in $W(C_0)^{epi}$, respectively, in $W(C_0),$
and we obtain a natural surjection
$W_P(E,C_0)\twoheadrightarrow B_P(E,C_0)$, $e\mapsto{\rm im}(\phi_e)$.

We put also
$$
\pi:=(pr_{2|P})^{-1}:\mathbb{P}^1\overset{\sim}\to P,\ \
\pi':=(pr_{2|P'})^{-1}:\mathbb{P}^1\overset{\sim}\to P',\ \
\rho:=(pr_{1|S})^{-1}:\mathbb{P}^1\overset{\sim}\to S.
$$

\begin{theorem}\label{phi}
Let $y_0\in l_0$ and $w\in\psi^{-1}(y_0)$ be  fixed points. Then,
for any vector bundle $E$ on $G$ there exists a morphism $f:\ X\to G$ such that:

(i) $f^*Q\simeq\mathcal{F}$, where $\FF$ is defined in \refth{F};

(ii) $f\circ\pi'=\psi$;

(iii) $f\circ\pi:\mathbb{P}^1\hookrightarrow G$ is an embedding such that
$(f\circ\pi)^*Q\simeq\mathcal{O}_{\mathbb{P}^1}(a_1)
\oplus...\oplus\mathcal{O}_{\mathbb{P}^1}(a_r),$
respectively,
$f\circ\rho:\mathbb{P}^1\hookrightarrow G$ is an embedding such that
$(f\circ\rho)^*Q\simeq (\mathcal{O}_{\mathbb{P}^1}(2))^r$;

(iv) $\mathbf{D}(E_{|f(P)})\le\mathbf{D}(E_{|C_0})$
and
$\mathbf{D}(E_{|f(S)})\le2r\mathbf{D}(E)$.
\end{theorem}

\unitlength=1.20mm \special{em:linewidth 0.4pt}
\linethickness{1.5pt}
\begin{picture}(128.00,157.00)
\bezier{168}(35.00,16.00)(55.00,22.00)(75.00,16.00)
\bezier{168}(75.00,16.00)(95.00,10.00)(115.00,16.00)
\bezier{140}(35.00,3.00)(47.00,8.00)(42.00,30.00)
\bezier{116}(42.00,30.00)(37.00,48.00)(45.00,55.00)
\bezier{212}(40.00,148.00)(40.00,118.00)(40.00,95.00)
\bezier{212}(110.00,148.00)(110.00,124.00)(110.00,95.00)
\bezier{200}(110.00,55.00)(110.00,30.00)(110.00,5.00)
\bezier{360}(35.00,100.00)(75.00,100.00)(125.00,100.00)
\linethickness{0.4pt} \put(110.00,136.00){\circle*{1.20}}
\put(110.00,124.00){\circle*{1.20}}
\put(110.00,112.00){\circle*{1.20}}
\put(8.00,100.00){\circle*{1.20}}
\put(27.00,124.00){\vector(-1,0){10.00}}
\bezier{120}(9.00,140.00)(24.00,142.00)(39.00,140.00)
\put(38.00,140.10){\vector(1,0){1.00}}
\bezier{400}(9.00,148.00)(59.00,152.00)(109.00,148.00)
\put(108.00,148.00){\vector(1,0){1.00}}
\put(110.00,100.00){\circle*{1.20}}
\put(110.00,81.00){\circle*{1.20}}
\put(40.00,81.00){\circle*{1.20}}
\put(122.00,84.00){\vector(0,1){13.00}}
\put(74.00,75.00){\vector(0,-1){12.00}}
\put(110.20,14.70){\circle*{1.20}}
\bezier{168}(35.00,43.00)(55.00,49.00)(75.00,43.00)
\bezier{168}(75.00,43.00)(95.00,37.00)(115.00,43.00)
\put(110.20,41.80){\circle*{1.20}}
\bezier{176}(115.00,45.00)(104.00,35.00)(75.00,40.00)
\bezier{176}(75.00,40.00)(52.00,46.00)(35.00,35.00)
\bezier{200}(115.00,48.00)(108.00,33.00)(75.00,36.00)
\bezier{184}(75.00,36.00)(47.00,40.00)(35.00,27.00)
\bezier{148}(111.00,15.00)(99.00,9.00)(75.00,12.00)
\bezier{168}(75.00,12.00)(49.00,17.00)(35.00,11.00)
\bezier{152}(111.00,15.00)(101.00,7.00)(75.00,9.00)
\bezier{172}(75.00,9.00)(45.00,14.00)(35.00,7.00)
\put(74.00,157.00){\makebox(0,0)[cc]{$X={\mathbb
P}^1\times{\mathbb P}^1$}}
\put(56.00,152.00){\makebox(0,0)[cc]{$\pi'$}}
\put(24.00,142.00){\makebox(0,0)[cc]{$\pi$}}
\put(22.00,125.00){\makebox(0,0)[cb]{$pr_2$}}
\put(56.00,137.00){\makebox(0,0)[cb]{$S_d$}}
\put(56.00,125.00){\makebox(0,0)[cb]{$S_2$}}
\put(56.00,113.00){\makebox(0,0)[cb]{$S_1$}}
\put(56.00,101.00){\makebox(0,0)[cb]{$S$}}
\put(45.00,131.00){\makebox(0,0)[cc]{$\vdots$}}
\put(114.00,137.00){\makebox(0,0)[cb]{$z_d$}}
\put(114.00,125.00){\makebox(0,0)[cb]{$z_2$}}
\put(114.00,113.00){\makebox(0,0)[cb]{$z_1$}}
\put(114.00,101.00){\makebox(0,0)[cb]{$ \ \ \ \ \ \ \ \ \ z_0=\pi'(w)$}}
\put(7.00,100.00){\makebox(0,0)[rc]{$w$}}
\put(8.00,93.00){\makebox(0,0)[cc]{${\mathbb P}^1$}}
\put(41.00,96.00){\makebox(0,0)[lc]{$P$}}
\put(106.20,96.00){\makebox(0,0)[lc]{$P'$}}
\put(74.00,96.00){\vector(0,-1){11.00}}
\put(75.00,90.00){\makebox(0,0)[lc]{$pr_1$}}
\put(123.00,90.00){\makebox(0,0)[lc]{$\rho$}}
\put(128.00,81.00){\makebox(0,0)[cc]{${\mathbb P}^1$}}
\put(40.00,79.00){\makebox(0,0)[cc]{$0$}}
\put(110.00,79.00){\makebox(0,0)[cc]{$\infty$}}
\put(75.00,69.00){\makebox(0,0)[lc]{$f$}}
\put(74.00,53.00){\makebox(0,0)[cc]{$f(X)$}}
\put(111.00,54.00){\makebox(0,0)[lc]}
\put(42.00,54.00){\makebox(0,0)[rc]}
\put(34.00,43.00){\makebox(0,0)[rc]{$f(S_d)$}}
\put(111.00,40.00){\makebox(0,0)[lc]{$f(z_1)=...=f(z_d)$}}
\put(111.00,13.00){\makebox(0,0)[lc]{$y_0=f(z_0)$}}
\put(34.00,34.00){\makebox(0,0)[rc]{$f(S_2)$}}
\put(34.00,26.00){\makebox(0,0)[rc]{$f(S_1)$}}
\put(44.00,43.00){\makebox(0,0)[cc]{$\vdots$}}
\put(59.00,21.00){\makebox(0,0)[cc]{$C^{2r}=f(S)$}}
\put(45.00,12.70){\makebox(0,0)[cc]{$.$}}
\put(45.00,12.10){\makebox(0,0)[cc]{$.$}}
\put(45.00,11.50){\makebox(0,0)[cc]{$.$}}
\put(36.00,1.00){\makebox(0,0)[cc]{$C^d=f(P)$}}
\put(110.00,1.00){\makebox(0,0)[cc]{$l_0=f(P')$}}
\bezier{228}(8.00,95.00)(8.00,132.00)(8.00,152.00)
\bezier{320}(35.00,136.00)(95.00,136.00)(115.00,136.00)
\bezier{320}(35.00,124.00)(105.00,124.00)(115.00,124.00)
\bezier{320}(35.00,112.00)(105.00,112.00)(115.00,112.00)
\bezier{396}(26.00,81.00)(115.00,81.00)(125.00,81.00)
\bezier{8}(10.00,100.00)(11.00,100.00)(12.00,100.00)
\bezier{8}(14.00,100.00)(15.00,100.00)(16.00,100.00)
\bezier{8}(18.00,100.00)(19.00,100.00)(20.00,100.00)
\bezier{8}(22.00,100.00)(23.00,100.00)(24.00,100.00)
\bezier{8}(26.00,100.00)(27.00,100.00)(28.00,100.00)
\bezier{8}(30.00,100.00)(31.00,100.00)(32.00,100.00)
\bezier{8}(40.00,83.00)(40.00,84.00)(40.00,85.00)
\bezier{8}(40.00,87.00)(40.00,88.00)(40.00,89.00)
\bezier{8}(40.00,91.00)(40.00,92.00)(40.00,93.00)
\bezier{8}(110.00,83.00)(110.00,84.00)(110.00,85.00)
\bezier{8}(110.00,87.00)(110.00,88.00)(110.00,89.00)
\bezier{8}(110.00,91.00)(110.00,92.00)(110.00,93.00)
\end{picture}

\vspace{0.5cm}

\begin{proof}
Recall (\refth{F}) that
\begin{equation}\label{F|P,P',S}
\mathcal{F}_{|P}\simeq\mathcal{O}_{\mathbb{P}^1}(a_1)\oplus...\oplus
\mathcal{O}_{\mathbb{P}^1}(a_r),\ \
\mathcal{F}_{|P'}\simeq\mathcal{O}_{\mathbb{P}^1}(d)\oplus(\mathcal{O}_{\mathbb{P}^1})^{r-1},\ \
\mathcal{F}_{|S}\simeq (\mathcal{O}_{\mathbb{P}^1}(2))^r.
\end{equation}
Furthermore, let us introduce the following notation:
\begin{equation*}\label{W_i}
z_0:=P'\cap S,\ \ \ W_0:=\Hom(V^\vee,H^0(\mathcal{F})),\ \ \
W_{-1}:=\Hom(V^\vee,H^0(\mathcal{F}(-1,0))),\ \ \
\end{equation*}
\begin{equation}\label{W_(S,P,PS)}
W_P:=\Hom(V^\vee,H^0(\mathcal{F}_{|P})),\ \ \
W_{-1,P}:=\Hom(V^\vee,H^0(\mathcal{F}(-1,0)_{|P}))\simeq W_P,\ \ \
\end{equation}
$$
W_{0,PS}:=\Hom(V^\vee,H^0(\mathcal{F}_{|P\cup S})),\ \ \
W_{-1,PS}:=\Hom(V^\vee,H^0(\mathcal{F}(-1,0)_{|P\cup S})),
$$
\begin{equation}\label{W_-1,S}
W_{-1,S}:=\Hom(V^\vee,H^0(\mathcal{F}(-1,0)_{|S}))\simeq
\Hom(V^\vee,H^0((\mathcal{O}_{\mathbb{P}^1}(1))^r)),\ \ \
\end{equation}
\begin{equation}\label{W_P'}
W_{P'}:=\Hom(V^\vee,H^0(\mathcal{F}_{|P'}))\simeq
\Hom(V^\vee,H^0(\mathcal{O}_{\mathbb{P}^1}(d)\oplus(\mathcal{O}_{\mathbb{P}^1})^{r-1})).
\end{equation}
(The right-hand isomorphisms in (\ref{W_-1,S}) and (\ref{W_P'}) follow from
(\ref{F|P,P',S})).

By using (\ref{F|P,P',S}), we see that the functor $\mathcal{H}om_{\OO_X}(V^\vee\otimes\OO_X,-)$ applied to the exact sequence $0\to\mathcal{F}(-1,0)\to\mathcal{F}\to\mathcal{F}_{|P'}\to0$ yields and exact sequence
$0\to\mathcal{H}om_{\mathcal{O}_X}(V^\vee\otimes\mathcal{O}_X,\mathcal{F}(-1,0))
\to\mathcal{H}om_{\mathcal{O}_X}(V^\vee\otimes\mathcal{O}_X,\mathcal{F})\to
\mathcal{H}om_{\mathcal{O}_{P'}}(V^\vee\otimes\mathcal{O}_{\mathbb{P}^1},
\mathcal{O}_{\mathbb{P}^1}(d)\oplus(\mathcal{O}_{\mathbb{P}^1})^{r-1})\to0.$
Passing to cohomology and using (\ref{h^1F}), we obtain the exact sequence
\begin{equation}\label{2 res P'}
0\to W_{-1}\overset{i_W}\to W_0\overset{res_{P'}}\to W_{P'}\to0.
\end{equation}
Next, for any $s\in W_{-1}$ we consider the composition morphism
$e_s:\ V^\vee\otimes\mathcal{O}_X\overset{s\otimes id}\to
H^0(\mathcal{F}(-1,0))\otimes\mathcal{O}_X\overset{ev}\to\mathcal{F}(-1,0)$,
and, for any
$z\in X\smallsetminus P'$, we consider the composition
\begin{equation*}\label{e_s(z)}
e_s(z):\ V^\vee\otimes\mathcal{O}_X\overset{e_s}\to\mathcal{F}(-1,0)
\overset{res_z}\twoheadrightarrow\mathcal{F}(-1,0)\otimes\mathbf{k}(z)\simeq\mathbf{k}(z)^r.
\end{equation*}
Passing to sections in the exact sequence
$0\to\mathcal{F}(-2,-1)\to\mathcal{F}(-2,0)\oplus\mathcal{F}(-1,-1)\to\mathcal{F}(-1,0)
\overset{res_z}\to\mathcal{F}(-1,0)\otimes\mathbf{k}(z)\to0$
and using (\ref{h^1F}) and (\ref{h^2F}), we obtain an epimorphism
$H^0(\mathcal{F}(-1,0))\overset{res_z}\twoheadrightarrow
H^0(\mathcal{F}(-1,0)\otimes\mathbf{k}(z))\simeq\mathbf{k}^r,$
and hence an induced epimorphism
\begin{equation}\label{r(z)}
r(z):\ W_{-1}\twoheadrightarrow
\Hom(V^\vee,H^0(\mathcal{F}(-1,0)\otimes\mathbf{k}(z)))
\simeq\Hom(V^\vee,\mathbf{k}(z)^r)=:W_z.
\end{equation}
Put
$Y(z):=\{s\in W_{-1}\ |\ e_s(z)\ {\rm is\ not\ surjective}\}$
and
$Y_0(z):=\{u\in W_z|\ u:V^\vee\to\mathbf{k}(z)^r\ {\rm is\ not\ surjective}\}.$
By definition $Y(z)=r(z)^{-1}(Y_0(z))$, and one easily checks that
$\codim_{W_z}Y_0(z)=\dim V-r+1.$ Therefore the surjectivity of $r(z)$ yields
\begin{equation}\label{codimY(z)}
\codim_{W_{-1}}Y(z)=\codim_{W_z}Y_0(z)=\dim V-r+1.
\end{equation}
If
$Y:=\underset{z\in X\smallsetminus P'}\cup Y(z)$, (\ref{codimY(z)}) implies
\begin{equation}\label{1 codimY(z)}
\codim_{W_{-1}}Y\ge \dim V-r-1>0.
\end{equation}
(Note that $\dim V-r-1>0$ according to our assumption from \refsec{sec2}.)

For each $z\in X\backslash P'$ the exact sequence (\ref{2 res P'}) and the map (\ref{r(z)}) fit in the diagram
\begin{equation*}
\xymatrix{
0\to W_{-1}\ \ \ \ \ar@{->>}[d]^{r(z)}\ar[r]^{i_W}&
W_0\ar@{->>}[d]^{r(z)}\ar[r]^{res_{P'}}&
W_{P'}\to0\\
W_z\ar@{=}[r]&
\Hom(V^\vee,H^0(\mathcal{F}\otimes\mathbf{k}(z))),
}
\end{equation*}
the right vertical map $r(z)$ being the natural restriction map. This diagram together with
(\ref{W_P'}) and the inequality (\ref{1 codimY(z)}) shows that, for the element
$$
g_\psi\in\Hom(V^\vee,H^0(\mathcal{O}_{\mathbb{P}^1}(d)
\oplus(\mathcal{O}_{\mathbb{P}^1})^{r-1}))=W_{P'}
$$
given in (\ref{g_psi}) and for a generic element
\begin{equation*}\label{epsilon}
\varepsilon\in res_{P'}^{-1}(g_\psi)\simeq W_{-1},
\end{equation*}
the composition
$
\varepsilon(z):\ V^\vee\otimes\mathcal{O}_X\overset{\varepsilon\otimes id}\to
H^0(\mathcal{F})\otimes\mathcal{O}_X\overset{ev}\to\mathcal{F}\overset{res_z}\to
\mathcal{F}\otimes\mathbf{k}(z)
$
is an epimorphism for any $z\in X\smallsetminus P'$. Moreover, since we can consider $g_{\psi}$ as an epimorphism
$:\ V^\vee\otimes\mathcal{O}_X\twoheadrightarrow\mathcal{F}_{|P'}$, $\varepsilon(z)$ is also an epimorphism for any $z\in P'$. This means that
$\varepsilon(z)$ is an epimorphism for any $z\in X$, i.e. that the morphism $\varepsilon:V^\vee\otimes\mathcal{O}_X\to\mathcal{F}$ is an epimorphism.
By the universality property of $G$ this means that there exists a morphism
$f=f_{\varepsilon}:X\to G$
such that
$\varepsilon=f_{\varepsilon}^*\gamma$, $\mathcal{F}\simeq f_{\varepsilon}^*Q$,
where $\gamma:V^\vee\otimes\mathcal{O}_G\to Q$ is the canonical epimorphism.
This together with (\ref{F|P,P',S}) yields (iii). In addition, since $\psi$ is determined by $g_\psi$, the equality
$g_\psi=res_{P'}(\varepsilon)$ directly implies (ii).

Next, we apply the functor
$\mathcal{H}om_{\mathcal{O}_X}(V^\vee\otimes\mathcal{O}_X,-)$
to the commutative diagram
$$
\xymatrix{
&0\ar[d]&0\ar[d]&0\ar[d]&\\
0\ar[r]&\mathcal{F}(-2,-1)\ar[r]\ar[d]&\mathcal{F}(-1,-1)\ar[r]\ar[d]&
\mathcal{F}(-1,-1)_{|P'}\ar[r]\ar[d]&0\\
0\ar[r]&\mathcal{F}(-1,0)\ar[r]\ar[d]&\mathcal{F}\ar[r]\ar[d]&\mathcal{F}_{|P'}\ar[r]\ar[d]&0\\
0\ar[r]&\mathcal{F}(-1,0)_{|P\cup S}\ar[r]\ar[d]&\mathcal{F}_{|P\cup S}\ar[r]\ar[d]&
\mathcal{F}\otimes\mathbf{k}(z_0)\ar[r]\ar[d]&0\\
&0&0&0&.
}
$$
Using (\ref{h^1F}), (\ref{h^2F}), (\ref{W_(S,P,PS)})-(\ref{W_P'}), we obtain the commutative diagram
\begin{equation}\label{diag3}
\xymatrix{
0\ar[r]&W_{-1}\ar[r]\ar@{->>}[d]^{res_{PS}}&W_0\ar[r]\ar@{->>}[d]&W_{P'}\ar[r]\ar@{->>}[d]^{r(z_0)}&0\\
0\ar[r]&W_{-1,PS}\ar[r]&W_{0,PS}\ar[r]^{res_{z_0}}&W_{z_0}\ar[r]&0\quad .
}
\end{equation}
Moreover, setting $\varepsilon_{z_0}:=r(z_0)(g_\psi)$, we have
\begin{equation*}\label{res^-1}
res_{z_0}^{-1}(\varepsilon_{z_0})\simeq W_{-1,PS}.
\end{equation*}
Similarly to (\ref{diag3}), using Theorem \ref{F} and
(\ref{W_(S,P,PS)})-(\ref{W_-1,S}), we obtain the surjective restriction maps
\begin{equation}\label{restrns 2}
res_PW_P\overset{res_P}\twoheadleftarrow W_{-1,PS}\overset{res_S}\twoheadrightarrow W_{-1,S}
\simeq\Hom(V^\vee,H^0((\mathcal{O}_{\mathbb{P}^1}(1))^r)).
\end{equation}
Now (\ref{diag3}) and (\ref{restrns 2})
together with Corollary \ref{cor B_S} show that
$$U:=res_{PS}^{-1}((res_P)^{-1}(W_P(E,C_0))\cap (res_S)^{-1}(W^*(E,y_0)))$$ is a dense open
subset of $W_{-1}$. Hence, for a generic element $\varepsilon\in U$ the corresponding morphism $f=f_{\varepsilon}:X\to G$ satisfies the conditions
(a) $e_{\pi}:=res_P(res_{PS}(\varepsilon))\in W_P(E,C_0)$
and
(b) $e_{\rho}:=res_S(res_{PS}(\varepsilon))\in W^*(E,y_0)$. Here, by construction, we have $\phi_{e_{\pi}}=f\circ\pi$ and,
respectively, $\phi_{e_{\rho}}=f\circ\rho$.
Now (a) together with (\ref{B_P(E,C_0)}) and (\ref{W_P(E,C_0)}) means that
$f(P)={\rm im}f\circ\pi={\rm im}\phi_{e_{\pi}}\in B_P(E,C_0)$, i.e. that
$\mathbf{D}(E_{|f(P)})\le\mathbf{D}(E_{|C_0})$;
respectively, (b) together with (\ref{W_S^*(E,y_0)}) and Corollary
\ref{B_S(E,y_0)} means that
$\mathbf{D}(E_{|f(S)})\le2r\mathbf{D}(E)$.
This yields (iv). The claim (i) is clear from the construction.
\end{proof}

\vspace{1cm}

\section{Main result}
\label{Quadric}
\vspace{0.5cm}

We now proceed to the main construction. Consider a twisted ind-Grassmannian $\GG$ defined by \refeq{eq1}, together with a vector bundle on $\mathbf{G}$, i.e. with vector bundles $E_m$ on $G_m:=G(r_m,V^{n_m})$ of fixed rank such that $\phi_m^*E_{m+1}=E_m$ for $m\geq 1$. We assume that $G_m$ is not isomorphic to a
projective space for any $m$. Set
\begin{equation*}\label{Phi_i}
\Phi_m:=\phi_{m-1}\circ...\circ\phi_1: G_1\hookrightarrow G_m.
\end{equation*}
The basic assumption \refeq{tame} that $\mathbf{G}$ is sufficiently twisted can be rewritten as
\begin{equation}\label{basic assumptn}
\underset{m\to\infty}\lim\ \ \frac{r_m}{\deg\Phi_m}=0.
\end{equation}
Note that \refeq{basic assumptn} is always satisfied if $r_m$ doesn't depend on $m$.

For a given integer $m>1$, fix a line $l_0$ in $G_m$ such that
\begin{equation}\label{l_0-max}
\mathbf{D}(E_{m|l_0})=\underset{{\rm lines}\ l\subset G_m}\max\ \
\mathbf{D}(E_{m|l})=:\mathbf{D}_m.
\end{equation}
 Let also $l_1$ be any line in $G_1$ such that
\begin{equation}\label{l_1-max}
\mathbf{D}(E_{1|l_1})=\underset{{\rm lines}\ l\subset G_1}\max\ \
\mathbf{D}(E_{1|l})=:\mathbf{D}_1.
\end{equation}
(Clearly, such lines $l_0$ and $l_1$ exist by semicontinuity). 

Put $C_0:=\Phi_m(l_1), \ \ \ d:=\deg C_{0}=\deg\Phi_m$. Fix a point $y_0\in l_0$ and a curve $C\in B(E,y_0)$ (see Theorem \ref{B_S(E,y_0)}) and fix
a degree $d$ morphism $\psi:\ \mathbb{P}^1\to l_0\hookrightarrow G_m\ \ $. Consider the
surface $X=\mathbb{P}^1\times\mathbb{P}^1$ together with the distinguished fibers $P,P'$ of the
projection $pr_1:X\to\mathbb{P}^1$ and $S$ of the projection $pr_1:X\to\mathbb{P}^1$. Applying
Theorem \ref{phi} to this datum, we obtain a morphism $f:X\to G_m$ such that, for
$E_X:=f^*E_m,$ we have:

(i) the morphism
$\mathbb{P}^1\overset{(pr_{2|P'})^{-1}}{\underset{\sim}\longrightarrow}P'\overset{f}\to G_m$
coincides with $\psi$, hence by (\ref{l_0-max}) and (\ref{l_1-max}) the vector bundle $E_X$
satisfies the equality
\begin{equation}\label{dD_m}
\mathbf{D}({E_X}_{|P'})=\deg(f_{|P'})\mathbf{D}(E_{m|l_0})=(\deg\Phi_m)\mathbf{D}_m;
\end{equation}

(ii) $f_{|P}$ and $f_{|S}$ are embeddings such that
\begin{equation}\label{{E_X}_{|P})}
\mathbf{D}({E_X}_{|P})=\mathbf{D}(E_{m|f(P)})\le\mathbf{D}(E_{m|C_0})=
\mathbf{D}(E_{m|\Phi_m(l_1)})=\mathbf{D}(E_{1|l_1})=\mathbf{D}_1,
\end{equation}
\begin{equation}\label{d(E_X|S)}
\mathbf{D}({E_X}_{|S})=\mathbf{D}(E_{m|f(S)})\le2r_m\mathbf{D}_m.
\end{equation}
Now applying the inequality (3.11) from \cite{DP} to (\ref{dD_m})-(\ref{d(E_X|S)}) we obtain
\begin{equation*}\label{main ineqs}
(\deg\Phi_m)\mathbf{D}_m=d\mathbf{D}_m\le4\rk {E_X}(\mathbf{D}_1+2)(2r_m\mathbf{D}_m+1)-
2\rk {E_X}.
\end{equation*}
But this inequality clearly contradicts to (\ref{basic assumptn}) for large enough $m$ if
$\mathbf{D}_m\ne0$.
Hence $\mathbf{D}_m=0$. Therefore, by \cite[Proposition 4.1]{PT} $E_m$ is trivial.
We thus have proved our main result.

\begin{theorem}\label{main}
There are no nontrivial vector bundles of finite rank on a sufficiently twisted ind-Grassmannian
$\mathbf{G}$.
\end{theorem}

We conclude this paper by a class of natural examples of twisted ind-Grassmannians (\ref{eq1})
with $r_m=r=\const$. Recall that any embedding of $\mathbb{P}^n$ to $\mathbb{P}^{n'}$ is given
by a subsystem of a complete linear system, i.e. by a composition of a Veronese embedding of
$\mathbb{P}^n$ into $\mathbb{P}^{n''},\ n''\ge n',$ and subsequent projection of
$\mathbb{P}^{n''}$ to $\mathbb{P}^{n'}$. In fact, this procedure extends to Grassmannians of
$r$-dimensional subspaces, $r$ being fixed. More precisely,
for each $m\ge1$ fix an integer $k_m>1$ and construct the Grassmannians $G_m=G(r,V_m)$
and their successive embeddings
$\phi_m: G_m\hookrightarrow G_{m+1}$
inductively by the following procedure.
Consider the flag variety $\Gamma_m=Fl(1,r;V_m)$ together with the natural embedding
$\Gamma_m\hookrightarrow P(V_m)\times G_m$ and the sheaf
$\mathcal{O}_\Gamma(1,k_m):=
(\mathcal{O}_{P(V_m)}(1)\boxtimes\mathcal{O}_{G_m}(k_m))_{|\Gamma_m}$.
Set $W_{m+1}:=H^0(\mathcal{O}_{\Gamma_m}(1,k_m))^\vee$. The embedding
$\theta_m:\Gamma_m\hookrightarrow P(W_{m+1})$ by the complete linear series
$|\mathcal{O}_{\Gamma_m}(1,k_m)|$
is (by construction) induced by a homogeneous embedding (in the sense of \cite[Sect. 4]{DP})
$\psi_m:G_m\hookrightarrow G(r,W_{m+1})$.
By composing $\psi_m$ with a possible rational projection of the form
$\pi:G(r,W_{m+1})\dasharrow G_{m+1}=G(r,V_{m+1})$, where $V_{m+1}$ is an appropriate quotient
of $W_{m+1}$,
we obtain an embedding $\phi_m:G_m\hookrightarrow G_{m+1}$.

Note that Theorem \ref{main} was proved in \cite{DP} for ind-Grassmannians \refeq{eq1} defined via certain homogeneous morphisms $\phi_m$ called twisted extensions, see \cite[Sect. 4.2]{DP}. It is not difficult to check that the above contructed embeddings
$\phi_m$ are not twisted extensions, hence Theorem \ref{main} is new for the corresponding
ind-Grassmannians $\mathbf G$.

\end{document}